

\magnification=1100
\overfullrule0pt

\input prepictex
\input pictex
\input postpictex
\input amssym.def


\def\qed{\hbox{\hskip 1pt\vrule width4pt height 6pt depth1.5pt \hskip 1pt}}

\def\CC{{\Bbb C}}
\def\RR{{\Bbb R}}
\def\ZZ{{\Bbb Z}}

\def\cB{{\cal B}}

\def\cF{{\cal F}}
\def\cG{{\cal G}}
\def\cT{{\cal T}}
\def\fg{{\goth g}}
\def\Ad{{\rm Ad}}

\def\Lie{{\rm Lie}}

\def\Ind{{\rm Ind}}
\def \supp{{\rm supp}}

 
\font\smallcaps=cmcsc10
\font\titlefont=cmr10 scaled \magstep1

\font\tinyrm=cmr10 at 8pt


\newcount\sectno
\newcount\subsectno
\newcount\resultno

\def\section #1. #2\par{
\sectno=#1
\resultno=0
\bigskip\noindent\centerline{\smallcaps #1.  #2}~\medbreak}

\def\subsection #1\par{\medskip\noindent{\bf  #1} }


\def\prop{ \global\advance\resultno by 1
\bigskip\noindent{\bf Proposition \the\sectno.\the\resultno. }\sl}
\def\lemma{ \global\advance\resultno by 1
\bigskip\noindent{\bf Lemma \the\sectno.\the\resultno. }
\sl}
\def\remark{ \global\advance\resultno by 1
\bigskip\noindent{\bf Remark \the\sectno.\the\resultno. }}
\def\example{ \global\advance\resultno by 1
\bigskip\noindent{\bf Example \the\sectno.\the\resultno. }}
\def\cor{ \global\advance\resultno by 1
\bigskip\noindent{\bf Corollary \the\sectno.\the\resultno. }\sl}
\def\thm{ \global\advance\resultno by 1
\bigskip\noindent{\bf Theorem \the\sectno.\the\resultno. }\sl}
\def\defn{ \global\advance\resultno by 1
\bigskip\noindent{\it Definition \the\sectno.\the\resultno. }\slrm}
\def\endthm{\rm\bigskip}

\def\endprop{\rm\bigskip}

\def\pf{\rm\bigskip\noindent{\it Proof. }}
\def\endpf{\qed\hfil\bigskip}


\def\formula{\global\advance\resultno by 1
\eqno{(\the\sectno.\the\resultno)}}
\def\formulano{\global\advance\resultno by 1 (\the\sectno.\the\resultno)}
\def\tableno{\global\advance\resultno by 1
\the\sectno.\the\resultno. }
\def\lformula{\global\advance\resultno by 1
\leqno(\the\sectno.\the\resultno)}

\def\monthname {\ifcase\month\or January\or February\or March\or April\or
May\or June\or
July\or August\or September\or October\or November\or December\fi}

\newcount\mins  \newcount\hours  \hours=\time \mins=\time
\def\now{\divide\hours by60 \multiply\hours by60 \advance\mins by-\hours
     \divide\hours by60         
     \ifnum\hours>12 \advance\hours by-12
       \number\hours:\ifnum\mins<10 0\fi\number\mins\ P.M.\else
       \number\hours:\ifnum\mins<10 0\fi\number\mins\ A.M.\fi}


\nopagenumbers
\def\runningtitle{\smallcaps Affine Hecke algebras}
\headline={\ifnum\pageno>1\eoheadline\else\firstheadline\fi}
\def\names{\smallcaps arun ram}
\def\firstheadline{\noindent Preliminary Draft \hfill  \today}
\def\firstheadline{}
\def\eoheadline{\ifodd\pageno\oddheadline\else\evenheadline\fi}
\def\oddheadline{\tenrm\hfil\runningtitle\hfil\folio}
\def\evenheadline{\tenrm \folio\hfil{\names}\hfil}

\vphantom{$ $}  
\vskip.75truein

\centerline{\titlefont Representations of rank two affine Hecke algebras}
\bigskip
\centerline{\rm Arun Ram${}^\ast$}
\centerline{Department of Mathematics}
\centerline{Princeton University}
\centerline{Princeton, NJ 08544}
\centerline{{\tt rama@math.princeton.edu}}
\centerline{Preprint: August 5, 1998}

\footnote{}{\tinyrm ${}^\ast$ Research supported in part by National
Science Foundation grant DMS-9622985, and a Postdoctoral Fellowship
at Mathematical Sciences Research Institute.}

\bigskip



\bigskip\bigskip

This paper classifies and constructs explicitly all the irreducible
representations of affine Hecke algebras of rank two root systems.
The methods used to obtain this classification are primarily
combinatorial and are, for the most part, an application of the methods used in
[Ra1].   I have made special effort to describe how the classification here
relates to the classifications by Langlands parameters (coming from $p$-adic
group theory) and by indexing triples (coming from a $q$-analogue
of the Springer correspondence).    There are several reasons for doing the
details of this classification:
\item{(a)}  The proof of the one of the main results of [Ra1] depends on this
classification of representations for rank two affine Hecke algebras.
Specifically, in the proof of Proposition 4.4 of [Ra1], one needs to know exactly
which weights can occur in calibrated representations.   The reason that this
naturally depends on a rank two classification is outlined in (d) below.
\item{(b)}  The examples here illustrate (and clarify)
results of [Ra1], [KL], [CG], [BM], [Ev], [Kr], [HO1-2].   Much of the power of
the combinatorial methods which are now available is evident from the calculations
in this paper, especially when one compares with the effort needed in other sources
(for example [Xi], Chapt. 11).
\item{(c)}  The explicit information here can be very useful for obtaining
results on representations of $p$-adic groups (see, for example, [Lu3]).
\item{(d)}   One hopes that eventually there will be a combinatorial construction
of all irreducible representations of all affine Hecke algebras.  I expect that
such a construction will depend heavily on the rank two cases.  This idea is
analogous to the way that the rank two cases are the basic building blocks in the
presentations of Coxeter groups by ``braid'' relations and the presentations of
Kac-Moody Lie algebras (and quantum groups) by Serre relations. 

\smallskip\noindent
The first section of this paper is a review of definitions and basic results about
affine Hecke algebras and their representations.  A few additional lemmas are
proved in order to aid the proofs and constructions in later sections.  The
remainder of the sections detail the classification and construction of the
irreducible representations of affine Hecke algebras of types $A_1$, $A_1\times
A_1$, $A_2$, $C_2$ and $G_2$.  In each case I have indicated how the results
here relate to the ``Langlands classification'', the classification of
Kazhdan and Lusztig [KL], and the results in [Ra1].

\subsection{Acknowledgements.}

This paper is part of a series [Ra1-3] [RR1-2] on
representations of affine Hecke algebras.  During this work I have benefited from
conversations with many people. To choose only a few, there were discussions with
S. Fomin, F. Knop, L. Solomon, M. Vazirani and N. Wallach which played an
important role in my progress.  There were several times when
I tapped into J. Stembridge's fountain of useful knowledge about root systems. 
D.-N. Verma helped at a crucial juncture by suggesting that I look at the paper
of Steinberg.  G. Benkart was a very patient listener on many occasions.  H.
Barcelo, P. Deligne, T. Halverson, R. Macpherson and R. Simion all gave large
amounts of time to let me tell them my story and every one of these sessions was
helpful to me in solidifying my understanding.

I single out Jacqui Ramagge with special thanks for everything she
has done to help with this project: from the most mundane typing and picture
drawing to deep intense mathematical conversations which helped to sort out
many pieces of this theory.  Her immense contribution is evident in
that some of the papers in this series on representations of affine Hecke
algebras are joint papers.  

A portion of this research was done during a semester long stay at Mathematical
Sciences Research Institute where I was supported by a Postdoctoral
Fellowship.  I thank MSRI and National Science Foundation for support of
my research.

\section 1.  Definitions and preliminary results

\subsection{The Weyl group.}

Let $R$ be a reduced irreducible root system in $\RR^n$, fix a set of positive
roots $R^+$ and let $\{\alpha_1,\ldots, \alpha_n\}$ be the corresponding
simple roots in $R$.  Let $W$ be the Weyl group corresponding to $R$.
Let $s_i$ denote the simple reflection in $W$ corresponding to the simple
root $\alpha_i$ and recall that $W$ can be presented by 
generators $s_1,s_2,\ldots, s_n$ and relations
$$\matrix{
s_i^2&=&1, &&\hbox{for $1\le i\le n$,} \cr
\underbrace{s_is_js_i\cdots}_{m_{ij} {\rm \ \ factors}}
&=& \underbrace{s_js_is_j\cdots}_{m_{ij} {\rm \ \ factors}}\;,
&\qquad &\hbox{for $i\ne j$,}  \cr
}
$$
where $m_{ij} = 
\langle \alpha_i,\alpha_j^\vee\rangle \langle\alpha_j,\alpha_i^\vee\rangle$,
and $\alpha_i^\vee =2\alpha_i/\langle \alpha_i,\alpha_i\rangle$.

\subsection{The Iwahori-Hecke algebra.}

Fix $q\in \CC^*$ such that $q$ is not a root of unity.  
The {\it Iwahori-Hecke algebra}
$H$ is the associative algebra over $\CC$ defined by 
generators $T_1,T_2,\ldots, T_n$ and relations
$$\matrix{
T_i^2&=& (q-q^{-1})T_i+1, &&\hbox{for $1\le i\le n$,} \cr
\underbrace{T_iT_jT_i\cdots}_{m_{ij} {\rm \ \ factors}}
&=& \underbrace{T_jT_iT_j\cdots}_{m_{ij} {\rm \ \ factors}}\;,
&\qquad &\hbox{for $i\ne j$,}  \cr
}
\formula$$
where $m_{ij}$ are the same as in the presentation of $W$.
For $w\in W$ define $T_w=T_{i_1}\cdots T_{i_p}$ where
$s_{i_1}\cdots s_{i_p} = w$ is a reduced
expression for $w$.  By [Bou, Ch. IV \S 2 Ex. 23], the element
$T_w$ does not depend on the choice of the reduced expression.
The algebra $H$ has dimension $|W|$ and the set $\{T_w\}_{w\in W}$
is a basis of $H$.

\subsection{The group $X$.}

The {\it fundamental weights} are the elements $\omega_1,\ldots, \omega_n$
of $\RR^n$ given by $\langle \omega_i, \alpha_j^\vee\rangle = \delta_{ij}$.
The {\it weight lattice}  is the $W$-invariant
lattice in $\RR^n$ given by
$$P= \sum_{i=1}^n \ZZ\omega_i.$$
Let $X$ be the
abelian group $P$ except written multiplicatively.  In other words,
$$X = \{ X^\lambda \ |\ \lambda\in P\},
\quad\hbox{and}\quad
X^\lambda X^\mu = X^{\lambda+\mu} = X^\mu X^\lambda,
\quad\hbox{for $\lambda,\mu\in P$}.$$
Let $\CC[X]$ denote the group algebra of $X$.
There is a $W$-action on $X$ given by
$wX^\lambda = X^{w\lambda}$
for $w\in W$, $X^\lambda\in X$,
which we extend linearly to  a $W$-action on 
$\CC[X]$.

\subsection{The affine Hecke algebra.}

The {\it affine Hecke algebra} $\tilde H$ associated to $R$ and $P$ is the
algebra given by
$$\tilde H  
= \CC\hbox{-span} \{ T_w X^\lambda \ |\ w\in W, X^\lambda\in X\}$$ 
where the multiplication of the $T_w$ is as in the Iwahori-Hecke algebra
$H$, the multiplication of the $X^\lambda$ is as in $\CC[X]$ and
we impose the relation
$$X^\lambda T_i = T_i X^{s_i\lambda} + (q-q^{-1})
{X^\lambda-X^{s_i\lambda}\over 1-X^{-\alpha_i}}, 
\qquad\hbox{for $1\le i\le n$ and $X^\lambda\in X$.}
\formula$$
This formulation of the definition of $\tilde H$ is due to Lusztig [Lu2]
following work of Bernstein and Zelevinsky.
The elements $T_wX^\lambda$, $w\in W$, $X^\lambda\in X$, form a basis
of $\tilde H$.

\subsection{Weights.}

Let 
$$T=\{ \hbox{group homomorphisms $t\colon X \to \CC^*$}\}.$$
The {\it torus} $T$ is an abelian group with a $W$-action given by
$(wt)(X^\lambda) = t(X^{w^{-1}\lambda})$.
For any element $t\in T$ define the {\it polar decomposition} 
$$t=t_rt_c,\qquad\hbox{ $t_r,t_c\in T$ such that $t_r(X^\lambda)\in \RR_{>0}$,
and $|t_c(X^\lambda)|=1$,}$$
for all $X^\lambda\in X$.  Let $Q^\vee = \sum_i \ZZ\alpha_i^\vee$.
There is a unique
$\mu\in
\RR^n$ and a unique 
$\nu\in \RR^n/Q^\vee$ such that
$$t_r(X^\lambda)=e^{\langle \mu,\lambda\rangle}
\qquad\hbox{and}\qquad
t_c(X^\lambda)=e^{2\pi i\langle \nu,\lambda\rangle},
\qquad\hbox{for all $\lambda\in P$.}\formula$$
In this way we identify the sets $T_r=\{ t\in T\ |\ t=t_r\}$
and $T_c=\{ t\in T\ |\ t=t_c\}$ with $\RR^n$
and $\RR^n/Q^\vee$, respectively.  
\subsection{Central characters.}

\thm (Bernstein, Zelevinsky, Lusztig [Lu1, 8.1])
The center of $\tilde H$ is  
$\CC[X]^W = \{ f\in \CC[X] \ |\ wf=f\}$.
\endthm

Since $\tilde H$ has countable dimension, Dixmier's version of Schur's lemma
implies that $Z(\tilde H)$ acts on an irreducible $\tilde H$-module $M$ by
scalars.   Let $t\in T$ be such that
$$pM = t(p)M, \qquad\hbox{for all $p\in Z(\tilde H)$.} $$
Since $Z(\tilde H)=\CC[X(T)]^W$ it follows that
$t(p(X))=(wt)(p(X))$ for all $w\in W$.  The $W$-orbit $Wt$ of
$t$ is the {\it central character} of $M$.  We shall often abuse notation
and refer to any weight $s\in Wt$ as ``the central character'' of $M$.

\subsection{Weight spaces.}

Let $M$ be a finite dimensional $\tilde H$-module.
For each $t\in T$ the {\it $t$-weight space of $M$} 
and the {\it generalized $t$-weight space} are the subspaces
$$\eqalign{
M_t &= \{ m\in M \ |\ X^\lambda m = t(X^\lambda)m
\hbox{\ for all $X^\lambda\in X$}\}
\qquad\hbox{and}  \cr
\cr
M_t^{\rm gen} &= \{ m\in M\ |\
\hbox{for each $X^\lambda\in X$,
$(X^\lambda-t(X^\lambda))^k m =0$ for some $k\in \ZZ_{>0}$}
\},  \cr
}
$$
respectively.
If $M_t^{\rm gen}\ne 0$ then $M_t\ne 0$.
In general $M\neq\bigoplus_{t\in T} M_t$, but we do have
$$
M=\bigoplus_{t\in T} M_t^{\rm gen}.
$$
This is a decomposition of $M$ into Jordan blocks
for the action of $\CC[X]$.  The set of {\it weights of $M$} is the set
$$\supp(M) = \{ t\in T \ |\ M_t^{\rm gen}\ne 0\}.\formula$$

\subsection{The calibration graph.}

Let $t\in T$.  Define a graph $\Gamma(t)$ with
$$
\matrix{
\hbox{Vertices:} \quad  &Wt,\hfil & \cr
\hbox{Edges:}  &wt \longleftrightarrow s_iwt,\quad
&\hbox{ if }\quad (wt)(X^{\alpha_i}) \ne q^{\pm2}. \cr
}
$$
\prop ([Ra1] Proposition 2.12) Let $M$ be a finite dimensional irreducible $\tilde
H$-module with central character $t$.  Then
$$\dim(M_s^{\rm gen})=\dim(M_{s'}^{\rm gen})$$
if $s$ and $s'$ are in the same connected component of the
calibration graph $\Gamma(t)$.
\endprop

If $t\in T$ define
$$P(t) = \{ \alpha>0\ |\ t(X^\alpha) = q^{\pm 2}\}
\qquad\hbox{and}\qquad
Z(t) = \{ \alpha>0\ |\ t(X^\alpha)=1 \}.\formula$$
For each subset $J\subseteq
P(t)$ define
$$\cF^{(t,J)} = 
\{ w\in W \ |\ R(w)\cap Z(t) = \emptyset,\ \ R(w)\cap P(t) = J\},
\formula$$ 
where $R(w) = \{ \alpha>0 \ |\ w\alpha<0 \}$ is the {\it inversion set} of $w$.
Define a {\it placed shape} to be a pair $(t,J)$ such that $t\in T$,
$J\subseteq P(t)$ and $\cF^{(t,J)}\ne \emptyset$.  The elements of the set
$\cF^{(t,J)}$ are called {\it standard tableaux} of shape $(t,J)$.

\prop ([Ra1] Theorem 2.14)
Let $t\in T$.  The connected components of the calibration graph $\Gamma(t)$ are
the sets 
$$\{ wt \ |\ w\in \cF^{(t,J)} \},
\qquad J\subseteq P(t),
\qquad\hbox{ such that $\cF^{(t,J)}\ne \emptyset$.}$$
\endprop

\subsection{Calibrated representations.}

A finite dimensional $\tilde H$-module  
$M$ is {\it calibrated} if $M_t^{\rm gen}=M_t$,
for all $t\in T$.

\prop ([Ra1] Proposition 4.2)
\item{(a)}  An irreducible $\tilde H$-module $M$ is calibrated if and only if
$\dim(M_t^{\rm gen}) = 1$ for all weights $t$ of $M$.
\item{(b)}  If $M$ is an irreducible $\tilde H$-module with regular central
character
$t$ (i.e. $Z(t)=\emptyset$) then $M$ is calibrated.
\endprop

Let $\alpha_i$ and $\alpha_j$ be simple roots in $R$ and let
$R_{ij}$ be the rank two root subsystem of $R$ which is generated
by $\alpha_i$ and $\alpha_j$.  Let $W_{ij}$ be the Weyl group of
$R_{ij}$, the subgroup of $W$ generated by the simple reflections
$s_i$ and $s_j$.  A weight $t\in T$ is {\it calibratable} for $R_{ij}$ if
one of the following two conditions holds:
\smallskip
\itemitem{(a)} $t(X^{\alpha})\ne 1$ for all $\alpha\in R_{ij}$,
\smallskip
\itemitem{(b)} $R_{ij}$ is of type $C_2$ or $G_2$ (assume that
$\alpha_i$ is the long root and $\alpha_j$ is the short root),
$ut(X^{\alpha_i})=q^2$ and $ut(X^{\alpha_j})=1$ for some $u\in W_{ij}$,
and $t(X^{\alpha_i})\ne 1$ and  $t(X^{\alpha_j})\ne 1$.
\smallskip\noindent
A {\it placed skew shape} is a placed shape $(t,J)$ such that 
for all $w\in \cF^{(t,J)}$ and all pairs $\alpha_i, \alpha_j$ of
simple roots in $R$ the weight $wt$ is calibratable for $R_{ij}$.

\thm  ([Ra1] Theorem 3.1 and Proposition 4.1)
\item{(a)}  Let $(t,J)$ be a placed skew shape and let $\cF^{(t,J)}$ be the
set of standard tableaux of shape $(t,J)$.  Define 
$$\tilde H^{(t,J)} =
\CC\hbox{-span} \{ v_w \ |\ w\in {\cF}^{(t,J)} \},$$
so that the symbols $v_w$ are a labeled basis of the
vector space $\tilde H^{(t,J)}$.
Then the following formulas make $\tilde H^{(t,J)}$
into an irreducible $\tilde H$-module:  For each $w\in {\cF}^{(t,J)}$,
$$
\matrix{
\hfill X^\lambda v_w &=& (wt)(X^\lambda) v_w, \hfill
&&\hbox{for $X^\lambda\in X$, and}\hfil \cr
\cr
\hfill T_i v_w 
&=& (T_i)_{ww} v_w + (q^{-1}+(T_i)_{ww}) v_{s_iw}, \hfill
&\qquad &\hbox{for $1\le i\le n$,}\hfil \cr
}
$$
where $\displaystyle{
(T_i)_{ww} = {q-q^{-1}\over 1 - (wt)(X^{-\alpha_i})}\;, }
$
and we set $v_{s_iw} = 0$ if $s_iw\not\in {\cF}^{(t,J)}$.
\item{(b)} If $M$ is an irreducible calibrated representation such that
$\supp(M)=\{ wt \ |\ w\in \cF^{(t,J)} \}$ for some placed skew shape
$(t,J)$ then $M$ is isomorphic to the module $\tilde H^{(t,J)}$ constructed in
(a).
\endthm

\remark
It follows from the results of Rodier [Ro] that if $M$ is an irreducible
$\tilde H$-module with regular central character (i.e. $Z(t)=\emptyset$) then $M$
satisfies the hypothesis of the statement of Theorem 1.11 (b). 
\endprop
 
\subsection{Langlands classification.}

The following discussion follows the work of Evens [Ev] and [KL, \S 8].
For this subsection it is convenient to assume that $q\in \RR_{>0}$ and $q\ne 1$.
For the general case see [KL, \S 8].  
Let $t\in T$ and let $t=t_rt_c$ be the polar decomposition of $t$.
Define 
$$\nu(t)\in \sum_{i=1}^n \RR\alpha_i^\vee
\qquad\hbox{by requiring}\qquad t_r(X^\lambda)= q^{2\langle
\lambda,\nu(t)\rangle},
\quad\hbox{for all $\lambda\in P$.}$$
A finite dimensional $\tilde H$-module $M$ is {\it tempered} if
for all weights $t$ of $M$ (as defined in (1.5)) we have
$$\langle\omega_i,\nu(t)\rangle\le 0,
\qquad\hbox{for all $1\le i\le n$.}$$
The module $M$ is {\it square integrable} if
$\langle \omega_i,\nu(t)\rangle<0$ for all $1\le i\le n$ and all weights $t$ of
$M$.


Let $I$ be a subset of the simple roots and let
$\tilde H_I$ be the subalgebra of $\tilde H$ generated by $T_i$, $i\in I$, and all
$X^\lambda\in X$.  We shall say that a finite dimensional
$\tilde H_I$-module is {\it tempered} if $I$ is the maximal set such that
for all weights $t$ of $M$,
$$\langle\omega_i,\nu(t)\rangle\le 0,
\qquad\hbox{for all $i\in I$.}$$

\thm  (see [Ev])  Let $I\subseteq \{1,2\ldots,n\}$ and let $\cT$ be an
irreducible tempered representation of $\tilde H_I$.

\smallskip\noindent
(a) $M_{\cT,I}=\Ind_{\tilde H_I}^{\tilde H} (\cT)$
has a unique irreducible quotient $L_{\cT,I}$.
\smallskip\noindent
(b) Every irreducible $\tilde H$-module is isomorphic to $L_{\cT,I}$
for some pair $(\cT,I)$.
\smallskip\noindent
(c) If $L_{\cT,I}\cong L_{\cT',I'}$ then $I=I'$ and $\cT\cong \cT'$ as
$\tilde H_I$-modules. 
\endthm

\noindent
The {\it Langlands parameters} of an irreducible $\tilde H$-module $M$
are given by the pair $(\cT,I)$ specified by Theorem 1.13 (b).

\subsection{Classification by indexing triples.}

Kazhdan and Lusztig [KL] (see also the important work of Ginzburg [CG])
gave a refinement of the Langlands classification.  Let $G$ be the simple
complex algebraic group with root system $R$ and weight lattice $P$.
An {\it indexing triple} $(s,u,\rho)$ consists of 
$$\matrix{ \hbox{a semisimple element $s\in G$,}\cr
\hbox{a unipotent element $u\in G$,} \cr
}
\qquad\hbox{such that}\qquad
sus^{-1} = u^{q^2}, $$
and an irreducible representation $\rho$ of
the component group
$A(s,u) = Z_G(s,u)/Z_G(s,u)^\circ$, where $Z_G(s,u)=Z_G(s)\cap Z_G(u)$.
Let $K(\cB_{s,u})$ be the $K$-theory of the variety
$$\cB_{s,u} = \{\hbox{Borel subgroups of $G$ containing both $s$ and $u$}\}.$$
By a theorem of Lusztig [Lu4] $K(\cB_{s,u})$ is an $\tilde H$-module.
The group $A(s,u)$ also acts on $K(\cB_{s,u})$ and this action commutes with
the action of $\tilde H$.  The {\it standard modules} $M_{s,u,\rho}$ are
the $\tilde H$-modules given by the decomposition 
$$K(\cB_{s,u}) = \bigoplus_{\rho} M_{s,u,\rho} \otimes \rho,$$
where the sum is over all irreducible representations of $A(s,u)$.

\thm [KL]
\item{(a)}  If $M_{s,u,\rho}\ne 0$ then it has a unique simple quotient
$L_{s,u,\rho}$.
\item{(b)}  Every simple $\tilde H$-module isomorphic to some
$L_{s,u,\rho}$.
\item{(c)} If $L_{s,u,\rho}\cong L_{s',u',\rho'}$ then there is a $g\in G$
such that $s'=gsg^{-1}$, $u'=gug^{-1}$, and $\rho'=\rho$.
\endthm

In this way each irreducible $\tilde H$-module corresponds to a unique 
(up to conjugation) indexing triple.  One can
replace  $u$ by $n = \ln u$ in the Lie algebra $\fg=\Lie(G)$ (see [CG, Ch. 8])
so that an indexing triple is 
$$\matrix{ \hbox{a semisimple element $s\in G$,}\cr
\hbox{a nilpotent element $n\in \fg$,} \cr
}
\qquad\hbox{such that}\qquad
\Ad(s) n = q^2 n, $$
and an irreducible representation $\rho$ of
the component group
$A(s,n) = Z_G(s,n)/Z_G(s,n)^\circ$, where $Z_G(s,n)=Z_G(s)\cap Z_G(n)$ and
$Z_G(n)$ is taken with respect to the adjoint action of $G$ on $\fg$.
We will use this form of the indexing triples in the examples in later sections.

\subsection{Principal series modules.}

Let $t\in T$ and let $\CC v_t$ be the one dimensional
$\CC[X]$-module corresponding to the character 
$t\colon X\to \CC^*$. 
Specifically, $\CC v_t$ is the one dimensional vector space with
basis $\{v_t\}$ and $\CC[X]$-action given by 
$$X^\lambda v_t =t(X^\lambda)v_t,
\qquad\hbox{for all $X^\lambda\in X$.}
$$  
The {\it principal series module} corresponding to $t$ is the
$\tilde H$-module
$$M(t) = \Ind_{\CC X}^{\tilde H}(\CC v_t).$$

\vskip-.1in

\thm [Ma]
\item{(a)}  Every irreducible $\tilde H$-module $M$ with central character $t$
is a composition factor of the principal series module $M(t)$.
\item{(b)}  If $w\in W$ and $t\in T$ then $M(t)$ and $M(wt)$ have the same
composition factors.
\endthm

\vskip-.15in

\thm (Kato's irreducibility criterion [Ka])  
Let $t\in T$ and let $P(t) = \{\alpha>0\ |\ t(X^\alpha) = q^{\pm 2}\}$.
The principal series module $M(t)$ is irreducible if and only if $P(t)=\emptyset$. 
\endthm

\noindent
{\bf Remark.}  Kato actually proves a more general result and thus needs
a further condition for irreducibility.  We have simplified matters by specifying
the weight lattice $P$ in our construction of the affine Hecke algebra.  One can
use any $W$-invariant lattice in $\RR^n$ and Kato works in this more general
situation.  When the one uses the weight lattice $P$, a result of Steinberg [St,
4.2, 5.3] says that the stabilizer $W_t$ of a point $t\in T$ under the action of
$W$ is always a reflection group. Because of this Kato's criterion takes a simpler
form.

\subsection{Weights of induced modules.} 

If $I\subseteq \{1,\ldots,n\}$ define $\tilde H_I$
to be the subalgebra of $\tilde H$ generated by $T_i$, $i\in I$,
and all $X^\lambda\in X$.  

\lemma  Let $t\in T$ such that $t(X^{\alpha_i})=q^2$ for all $i\in I$ 
and let $\CC v_t$ be the one dimensional
$\tilde H_I$-module with basis $\{v_t\}$ and $\tilde H_I$-action given by 
$$
T_iv_t = qv_t, \quad \hbox{for $i\in I$,\quad and}\qquad
X^\lambda v_t=t(X^\lambda)v_t,\quad \hbox{for all $X^\lambda\in X$.}
$$  
Let $W/W_I$ be the set of minimal length coset representatives of cosets of
$W_I$ in $W$.  Then the weights of the $\tilde H$-module 
$M=\Ind_{\tilde H_I}^{\tilde H}(\CC v_t)$ are $wt$,
$w\in W/W_I$, and
$$\dim(M_{wt}^{\rm gen})=\hbox{(\# of\ \  $u\in W/W_I$ such that $ut = wt$).}$$
\pf
The module $M$ has basis $\{T_w\otimes v_t\ |\ w\in W/W_I\}$.  By writing
$T_w=T_{i_1}\cdots T_{i_p}$ for a reduced word $w=s_{i_1}\ldots s_{i_p}$
and inductively using the defining relation (1.2) we get
$$\eqalign{
X^\lambda(T_w\otimes v_t) &= 
t(X^{w^{-1}\lambda})(T_w\otimes v_t) +\sum_{v<w} a_v(t)(T_v\otimes v_t) \cr
&=
(wt)(X^{\lambda})(T_w\otimes v_t) +\sum_{v<w} a_v(t)(T_v\otimes v_t), \cr}
$$
where the sum is over $v\in W$ which are less than $w$ in Bruhat order
and $a_v(t)\in \CC$.  This shows that the eigenvalues of $X^\lambda$ on $M$
are $(wt)(X^\lambda)$.  The result follows by counting the multiplicity of each
eigenvalue.
\endpf

\subsection{The $\tau$ operators.}

The maps $\tau_i\colon M_t^{\rm gen}\to M_{s_it}^{\rm gen}$ 
defined below are local operators on $M$ in the sense that they
act on each weight space $M_t^{\rm gen}$ of $M$ separately.  The operator
$\tau_i$ is only defined on weight spaces $M_t^{\rm gen}$ such that
$t(X^{\alpha_i})\ne 1$.  

\prop ([Ra1] Proposition 2.7)
Let $t\in T$ such that $t(X^{\alpha_i})\ne 1$ and let $M$
be a finite dimensional $\tilde H$-module.  Define
$$
\matrix{
\tau_i \colon 
&M_t^{\rm gen} &\longrightarrow &M_{s_it}^{\rm gen} \cr
\cr
& m & \longmapsto &
\displaystyle{
\left(T_i - {q-q^{-1}\over 1-X^{-\alpha_i} }\right) m.} \cr
}
$$
\item{(a)}
The map $\tau_i\colon M_t^{\rm gen} \longrightarrow M_{s_it}^{\rm gen}$
is well defined.
\smallskip
\item{(b)}  As operators on $M_t^{\rm gen}$, \quad
$\displaystyle{ X^\lambda \tau_i = \tau_i X^{s_i\lambda}}$,
for all $X^\lambda\in X$.
\smallskip
\item{(c)}  As operators on $M_t^{\rm gen}$, \quad
$\displaystyle{
\tau_i\tau_i=
{(q-q^{-1}X^{\alpha_i})(q-q^{-1}X^{-\alpha_i})\over 
(1-X^{\alpha_i})(1-X^{-\alpha_i}) }.
}
$
\smallskip
\item{(d)}  Let $1\le i\ne j\le n$ and let $m_{ij}$ be as in (1.1).
Then \quad   
$\displaystyle{\underbrace{\tau_i\tau_j\tau_i\cdots}_{m_{ij} {\rm\ factors}}
= \underbrace{\tau_j\tau_i\tau_i\cdots}_{m_{ij} {\rm\ factors}},
}$
\hfil\break
whenever both sides are well defined operators on $M_t^{\rm gen}$.
\endprop

\lemma Let $t\in T$ such that $t(X^{\alpha_i})=1$ and suppose that
$M$ is an $\tilde H$-module such that $M_t^{\rm gen}\ne 0$.  Let $W_t$ be the 
stabilizer of $t$ under the action of $W$ on $T$. 
Assume that $\bar w\in W/W_t$ is such that $t$ and $\bar wt$ are in the same
connected component of $\Gamma(t)$.   Let $w$ be a minimal length coset
representative of $\bar w$. Then 
\smallskip
\itemitem{(a)} $\dim(M_{wt}^{\rm gen})\ge 2$, and 
\smallskip
\itemitem{(b)} If $M_{s_jwt}^{\rm gen}=0$ then 
$(\bar wt)(X^{\alpha_j})=q^{\pm2}$ and $\langle
w^{-1}\alpha_j,\alpha_i^\vee\rangle =0$.
\pf
Let $M(t)$ be the two dimensional principal series module for the
affine Hecke algebra $\tilde HA_1$ of type $A_1$ (see \S 2 central character
$t_o$).   Then $M(t)=M(t)_t^{\rm gen}$ and has basis
$\{v_t, T_1v_t\}$. Let $n_t$ be a nonzero weight vector in
$M_t$.  There is a unique $\tilde HA_1$-module homomorphism
$$\matrix{
M(t) &\longrightarrow M \cr
v_t &\longmapsto n_t \cr
}$$
where we view $M$ as an $\tilde HA_1$-module by restriction to the parabolic
subalgebra $\tilde H_{\{i\}}\subseteq \tilde H$.  This homomorphism must
be an injection since $M(t)$ is irreducible.  Thus the vectors
$n_t, T_i n_t$ span a two dimensional subspace of
$M_t^{\rm gen}$ and $X^\lambda\in X$ acts on this subspace by the matrix
$$\phi_t(X^\lambda)
=t(X^\lambda)\pmatrix{ 1 &(q-q^{-1})\langle \lambda,\alpha_i^\vee\rangle \cr
0 &1\cr}.$$

Let $w=s_{i_1}\cdots s_{i_p}$ be a reduced expression of $w$.  Since $t$ and $wt$
are in the same connected component of $\Gamma(t)$ we can use Proposition 1.18 (c)
to show that the map
$$\tau_w=\tau_{i_1}\cdots \tau_{i_p}\colon M_t^{\rm gen}
\longrightarrow M_{wt}^{\rm gen}$$
is well defined and bijective.  Thus the vectors $\tau_w n_t,\tau_w T_i n_t$
span a two dimensional subspace of $M_{wt}^{\rm gen}$ and, by Proposition 1.18 (b) 
$X^\lambda\in X$ acts on this subspace by the matrix
$$\phi_{wt}(X^\lambda)
=t(X^{w^{-1}\lambda})\pmatrix{ 1 &(q-q^{-1})\langle
w^{-1}\lambda,\alpha_i^\vee\rangle
\cr 0 &1\cr}.$$
This proves (a).
Then
$$\phi_{wt}(1-X^{-\alpha_j})
=(1-t(X^{-w^{-1}\alpha_j}))
\pmatrix{ 1 
&\displaystyle{
{(q-q^{-1})t(X^{-w^{-1}\alpha_j}) \over 1-t(X^{-w^{-1}\alpha_j}) }
\langle -w^{-1}\alpha_j,\alpha_i^\vee\rangle  }  \cr 
0 &1\cr}.
$$
Since $M_{s_jwt}^{\rm gen} = 0$, 
$\tau_j\colon M_{wt}^{\rm gen} \to M_{s_jwt}^{\rm gen}$ is the zero map and 
so
$$\phi_{wt}(T_j) = \phi_{wt}\left({q-q^{-1}\over 1-X^{-\alpha_j}}\right)
={q-q^{-1}\over 1-t(X^{-w^{-1}\alpha_j})}
\pmatrix{ 1 
&\displaystyle{ {(q-q^{-1})t(X^{-w^{-1}\alpha_j})\over 1-t(X^{-w^{-1}\alpha_j})}
\langle w^{-1}\alpha_j,\alpha_i^\vee\rangle } \cr 
0 &1\cr}.
$$
The relation $T_j^2 = (q-q^{-1})T_j+1$ is the same as
$(T_j-q)(T_j+q^{-1}) = 0$.  This relation forces $\phi_{wt}(T_j)$ to
have Jordan blocks of size $1$ and eigenvalues $\pm q^{\pm1}$.  It follows
that $t(X^{w^{-1}\alpha_j}) = q^{\pm 2}$ and 
$\langle w^{-1}\alpha_j,\alpha_i^\vee\rangle = 0$.
\endpf

\section 2. Classification for $A_1$

\bigskip

The root system $R$ for $A_1$ has one simple root $\alpha_1$ and
fundamental weight $\omega_1={1\over2}\alpha_1$.

\bigskip
\subsection{Irreducible representations.}

Table 2.1 lists the irreducible $\tilde H$-modules by their
central characters.  The sets $P(t)$ and $Z(t)$ are as given in (1.7) and
correspond to the choice of representative for the central character displayed in
Figure 2.1.  The Langlands parameters usually consist of a pair $(\cT,I)$
where $I$ is a subset of $\{1\}$ and $\cT$ is a
tempered representation for the parabolic subalgebra $\tilde H_I$.  In our cases
the tempered representation $\cT$ of $\tilde H_I$ is completely determined by a
character $t\in T$.  Specifically, $\cT$ is the only tempered representation of
$\tilde H_I$ which has $t$ as a weight. In the labeling in
Table 2.1 we have replaced the representation $\cT$ by the weight $t$.  The
nilpotent element $e_{\alpha_1}$ is a representative of the root
space $\fg_{\alpha_1}$ for the Lie algebra $\fg={\goth{sl}}_2$.  For each
calibrated module with central character $t$ we have listed the subset $J\subseteq
P(t)$ such that $(t,J)$ is the corresponding placed skew shape (see Theorem
1.11).   The abbreviation `nc' indicates modules that are not calibrated. 
\bigskip
$$\matrix{
\hbox{Central} &P(t) &Z(t) &\hbox{Dimension}
&\hbox{Langlands} &\hbox{Indexing} &\hbox{Calibration} \cr
\hbox{character} &&&&\hbox{parameters} &\hbox{triple} &\hbox{set $J$} \cr
\cr
t_a &\{\alpha_1\} &\emptyset   
  &1  &(t_a,\emptyset) &(t_a,0,1) &\emptyset\cr
&&&1 &\hbox{tempered} &(t_a,e_{\alpha_1},1) &\{\alpha_1\} \cr
\cr
t_b &\emptyset &\emptyset   
  &2 &(t_b,\emptyset) &(t_b,0,1)  &\emptyset \cr
\cr
t_o &\emptyset &\{\alpha_1\}   
  &2 &\hbox{tempered} &(t_o,0,1) &\hbox{nc} \cr
}
$$
\medskip
\centerline{{\bf Table 2.1.}  Irreducible representations}

\bigskip

Figure 2.1 displays the real parts of the central
characters in Table 2.1.  If $t\in T$ then the polar decomposition $t=t_rt_c$
determines an element $\mu\in \RR^n$ such that
$t_r(X^\lambda)=e^{\langle\lambda,\mu\rangle}$ (see (1.3)).  For each central
character $t_p$ the point labeled by $p$ in Figure 2.1 is the graph of the
corresponding $\mu_p\in \RR^n$.  Assume (for pictorial convenience) that
$q$ is a positive real number and let
$$
H_{\alpha_1} = \{ x\in \RR \ |\ \langle \alpha_1,x\rangle = 0 \},
\qquad\hbox{and}\qquad
H_{\alpha_1\pm \delta} 
= \{ x\in \RR \ |\ \langle \alpha_1,x\rangle = \ln(q^{\pm2}) \}.
$$
The $\vert$ marks indicate the
(affine) hyperplanes $H_{\alpha_1\pm\delta}$.
$$
\beginpicture
\setcoordinatesystem units <0.75cm,0.75cm>         
\setplotarea x from -4 to 4, y from -0.5 to 0.5    
\linethickness=0.5pt                          
\putrule from -4 0 to 4 0          
\put{$+$} at 2 0
\put{$+$} at -2 0
\put{$\bullet$} at 2 0
\put{$t_a$}[b] at 2 .2
\put{$\bullet$} at 3 0
\put{$t_b$}[b] at 3 .2
\put{$\bullet$} at 0 0
\put{$t_o$}[b] at 0 .2
\put{$H_{\alpha_1}$}[t] at 0 -0.25
\put{$H_{\alpha_1+\delta}$}[t] at 2 -0.25
\put{$H_{\alpha_1-\delta}$}[t] at -2 -0.25
\endpicture
$$
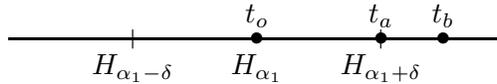
\centerline{{\bf Figure 2.1.}  Real parts of central characters in Table 2.1}

\medskip
\subsection{Tempered and square integrable representations.}

The tempered (resp. square integrable) $\tilde H$-modules are the ones which have
all their weight spaces in the closure (resp. interior) of the dotted region of
Figure 2.2. 
$$
\beginpicture
\setcoordinatesystem units <0.75cm,0.75cm>         
\setplotarea x from -4 to 4, y from -0.5 to 0.5    
\linethickness=0.5pt                          
\putrule from 4 0 to 0 0          
\put{$\bullet$} at 0 0
\put{$\bullet$} at -2 0
\put{$\vert$} at -2 0
\put{$\vert$} at 2 0
\put{$t_o$}[b] at 0 0.25
\put{$s_1t_a$}[b] at -2 0.25
\setdots
\plot 0 0 -4 0 /
\endpicture
$$
\centerline{{\bf Figure 2.2.} Real parts of weights of tempered representations}

\bigskip
The irreducible tempered representations with real central character can be indexed
by the irreducible representations of the symmetric group $S_2$
(see [BM]).  These representations are indexed by the partitions $(2),
(1^2)$ of $2$.   Let $e_{\alpha_1}$ be an element of the root space
$\fg_{\alpha_1}$ for the Lie algebra $\fg={\goth{sl}}_2$.  The two nilpotent
orbits in $\fg$ and the corresponding tempered representations of
$\tilde H$ are as in Table 2.2.  
$$\matrix{
\hbox{Nilpotent orbit}   &\hbox{Indexing triple}
&\hbox{Square integrable} &\hbox{$W$ representation}\cr
\hbox{regular}  &(t_a, e_{\alpha_1},1) &\hbox{yes}
&(1^2)\cr
0   &(t_o,0,1) &\hbox{no} &(2) \cr
}$$
\centerline{{\bf Table 2.2.}  Tempered representations and the Springer
correspondence}

\bigskip\noindent

$$
\beginpicture
\setcoordinatesystem units <0.75cm,0.75cm>         
\setplotarea x from -2 to 2, y from -0.5 to 1    
\linethickness=0.5pt                          
\put{$\bullet$} at .5 0
\put{$t_a$}[b] at .5 0.2
\put{$\bullet$} at -.5 0
\put{$s_1t_a$}[b] at -.5 0.2
\endpicture
\qquad\qquad
\beginpicture
\setcoordinatesystem units <0.75cm,0.75cm>         
\setplotarea x from -2 to 2, y from -0.5 to 1    
\linethickness=0.5pt                          
\put{$\bullet$} at .5 0
\put{$t_b$}[b] at .5 0.2
\put{$\bullet$} at -.5 0
\put{$s_1t_b$}[b] at -.5 0.2
\putrule from -.5 0 to .5 0          
\endpicture
\qquad\qquad
\beginpicture
\setcoordinatesystem units <0.75cm,0.75cm>         
\setplotarea x from -1 to 1, y from -0.5 to 1    
\linethickness=0.5pt                          
\put{$\bullet$} at 0 0
\put{$t_o$}[b] at 0 0.2
\endpicture
$$
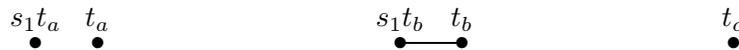
\centerline{{\bf Figure 2.3.}  Calibration graphs for central characters in
Table 2.1}

\medskip

\subsection{The analysis.}

\smallskip\noindent
{\it Central character $t_a$:}
There are two one-dimensional representations $\CC v_a$ and $\CC v_{s_1a}$ with
central character $t_a$.  These representations are given explicitly by
$$\eqalign{
X^\lambda v_a &= t_a(X^\lambda)v_a, \cr
T_1v_a &= qv_a, \cr}
\qquad\hbox{and}\qquad
\eqalign{
X^\lambda v_{s_1a} &= (s_1t_a)(X^\lambda)v_{s_1a}, \cr
T_1v_{s_1a}&= -q^{-1}v_{s_1a}, \cr
}$$
respectively.  One uses Theorem 1.15 and the fact that the principal series module
$M(t_a)$ is two dimensional to conclude that $\CC v_a$ and $\CC
v_{s_1a}$ are the only irreducible representations of $\tilde H$ with central
character $t_a$.
\smallskip\noindent
{\it Central character $t_b$:}  By Theorem 1.15 and Kato's irreducibility
criterion, Theorem 1.16, the only irreducible representation with central character
$t_b$ is the principal series module $M(t_b)$.  Alternatively, this module can be
constructed by applying Theorem 1.11 to the placed skew shape $(t_b,\emptyset)$.
\smallskip\noindent
{\it Central character $t_o$:}  The weights given by $t_o(X^{\omega_1})=\pm 1$
are the {\it two} central characters $t_o\in T$ which satisfy $P(t)=\emptyset$,
$Z(t)=\{\alpha_1\}$.  In either case Kato's irreducibility criterion (Theorem 1.16)
tells us that the principal series module $M(t_o)$ is irreducible.
This module has basis $\{v_t, T_1v_t\}$ and action given by
$$
\phi(X^\lambda)
=t(X^\lambda)\pmatrix{ 1 &(q-q^{-1})\langle \lambda,\alpha_1^\vee\rangle \cr
0 &1\cr} 
\qquad\hbox{and}\qquad
\phi(T_1) = 
\pmatrix{0 &1\cr 1 &q-q^{-1} \cr}.$$

\section 3. Classification for $A_1\times A_1$

The affine Hecke algebra of $A_1\times A_1$ is naturally isomorphic to 
$\tilde HA_1\otimes \tilde HA_1$.
The finite dimensional irreducible representations of $\tilde HA_1\otimes
\tilde HA_1$ are all of the form $M\otimes N$ where $M$ and $N$ are finite
dimensional irreducible representations of $\tilde HA_1$.

\vfill\eject

\section 4. Classification for $A_2$

\bigskip

The root system $R$ for $A_2$ has simple roots $\alpha_1$ and $\alpha_2$,
fundamental weights $\omega_1$ and $\omega_2$, and 
$$
\eqalign{
\langle \alpha_1,\alpha_2^\vee\rangle &= -1 \cr
\langle \alpha_2,\alpha_1^\vee\rangle &= -1, \cr}
\qquad\eqalign{
\omega_1 &= \hbox{$1\over3$}(2\alpha_1+\alpha_2)\cr
\omega_2&=\hbox{$1\over3$}(\alpha_1+2\alpha_2), \cr }
\qquad\hbox{and}\qquad
\eqalign{
\alpha_1 &= 2\omega_1-\omega_2 \cr 
\alpha_2&=-\omega_1+2\omega_2. \cr}
$$

\subsection{Irreducible representations.}

Table 4.1 lists the irreducible $\tilde H$-modules by their
central characters.  The sets $P(t)$ and $Z(t)$ are as given in (1.7) and
correspond to the choice of representative for the central character displayed in
Figure 4.1.  The Langlands parameters usually consist of a pair $(\cT,I)$
where $I$ is a subset of $\{1,2\}$ and $\cT$ is a
tempered representation for the parabolic subalgebra $\tilde H_I$.  In our cases
the tempered representation $\cT$ of $\tilde H_I$ is completely determined by a
character $t\in T$.  Specifically, $\cT$ is the only tempered representation of
$\tilde H_I$ which has $t$ as a weight. In the labeling in
Table 4.1 we have replaced the representation $\cT$ by the weight $t$.  The
nilpotent elements $e_{\alpha_1}$ and $e_{\alpha_2}$ are representatives of the
root spaces $\fg_{\alpha_1}$ and $\fg_{\alpha_2}$, respectively, where $\fg$ is
the Lie algebra $\fg={\goth{sl}}_3$.  For each calibrated module with central
character $t$ we have listed the subset
$J\subseteq P(t)$ such that $(t,J)$ is the corresponding placed skew shape (see
Theorem 1.11).   The abbreviation `nc' indicates modules that are not calibrated.
$$
\matrix{
\hbox{Central} &P(t) &Z(t) &\hbox{Dimension} &\hbox{Langlands} 
&\hbox{Indexing} &\hbox{Calibration} \cr
\hbox{character} & & & &\hbox{parameters} &\hbox{triple} &\hbox{set $J$}
\cr
\cr
t_a &\{\alpha_1,\alpha_2\} &\emptyset 
&1 &(t_a,\emptyset) &(t_a,0,1) &\emptyset\cr 
&&&2 &(s_1t_a,\{2\}) &(t_a,e_{\alpha_2},1) &\{\alpha_2\}\cr
&&&2 &(s_2t_a,\{1\}) &(t_a,e_{\alpha_1},1) &\{\alpha_1\}\cr
&&&1 &\hbox{tempered} &(t_a,e_{\alpha_1}+e_{\alpha_2},1) 
&\{\alpha_1,\alpha_2\}\cr
\cr
t_b &\{\alpha_2\} &\emptyset 
&3 &(t_b,\emptyset) 
&(t_b,0,1) &\emptyset\cr 
&&&3 &(s_2t_b,\{2\})^\dagger &(t_b,e_{\alpha_2},1) 
&\{\alpha_2\}\cr
\cr
t_c &\{\alpha_2,\alpha_1+\alpha_2\} &\{\alpha_1\}
&3 &(t_c,\{1\}) &(t_c,0,1) &\hbox{nc}\cr  
&&&3 &(s_2t_c,\{2\})  &(t_c,e_{\alpha_2},1) &\hbox{nc}\cr
\cr
t_d &\{\alpha_1,\alpha_1+\alpha_2\} &\{\alpha_2\}
&3 &(t_d,\{2\}) &(t_d,0,1) &\hbox{nc}\cr  
&&&3 &(s_1t_d,\{1\})  &(t_d,e_{\alpha_1},1) &\hbox{nc}\cr
\cr
t_e &\emptyset &\{\alpha_1\}
&6 &(t_e,\{1\}) &(t_e,0,1) &\hbox{nc} \cr
\cr
t_f &\emptyset &\{\alpha_2\}
&6 &(t_f,\{2\}) &(t_f,0,1) &\hbox{nc} \cr
\cr
t_g &\emptyset &\emptyset
&6 &(t_g,\emptyset) &(t_g,0,1) &\emptyset \cr
\cr
t_o &\emptyset &\{\alpha_1,\alpha_2\}
&6 &\hbox{tempered} &(t_o,0,1) &\hbox{nc} \cr
\cr
}
$$
\smallskip
\centerline{{\bf Table 4.1.}  Irreducible representations}
\footnote{}{$^\dagger$ There is one case when
this representation is tempered, see Table 4.2.}

\bigskip
Figure 4.1 displays the real parts of the central
characters in Table 4.1.  If $t\in T$ then the polar decomposition
$t=t_rt_c$ determines an element $\nu\in \RR^n$ such that
$t_r(X^\lambda) = e^{\langle \lambda, \nu \rangle}$ (see (1.3)).  
For each central character $t_p$ the point labeled by $p$ in Figure 4.1 is
the graph of the corresponding $\nu_p\in \RR^n$.  Assume (for pictorial
convenience) that $q$ is a positive real number and let
$$
H_\beta = \{ x\in \RR^n \ |\ \langle \beta,x\rangle = 0 \},
\qquad\hbox{and}\qquad
H_{\beta\pm \delta} = \{ x\in \RR^n \ |\ \langle \beta,x\rangle = \ln(q^{\pm2})
\},
$$
for each positive root $\beta$.  The dotted lines display the (affine) hyperplanes
$H_{\beta\pm\delta}$.
$$
\beginpicture
\setcoordinatesystem units <1cm,0.866cm>         
\setplotarea x from -4 to 4, y from -3 to 3.5    
\put{$H_{\alpha_1}$}[b] at 0 3.05
\put{$H_{\alpha_2}$}[l] at 2.7 1.6
\put{$H_{\alpha_1+\alpha_2}$}[l] at 2.7 -1.6
\put{$\bullet$} at .8660 1.5
\put{$a$}[bl] at .93 1.70
\put{$\bullet$} at .4330 1.25
\put{$b$}[b] at .4330 1.45
\put{$\bullet$} at 0 1
\put{$c$}[br] at -0.1 1.2
\put{$\bullet$} at .8660 0.5
\put{$d$}[l] at 1.1 0.55
\put{$\bullet$} at 0 2
\put{$e$}[r] at -0.1 2
\put{$\bullet$} at 1.7320 1
\put{$f$}[tl] at 1.75 .9
\put{$\bullet$} at 1.75 1.6
\put{$g$}[l] at 1.9 1.6
\put{$\bullet$} at 0 0
\put{$o$}[r] at -0.25 0
\plot 0 3    0 -3 /
\plot -2.5980 -1.5    2.5980 1.5 /
\plot -2.5980  1.5    2.5980 -1.5 /
\setdots
\plot  .8660 2.4340  .8660 -2.6340 /
\plot -.8660 2.4340 -.8660 -2.6340 /
\plot -2.5980 -0.5    1.7320 2.0 /
\plot -1.7320  2.0    2.5980 -0.5 /
\plot -1.7320  -2.0   2.5980 0.5 /
\plot -2.5980  0.5    1.7320 -2.0 /
\endpicture
$$
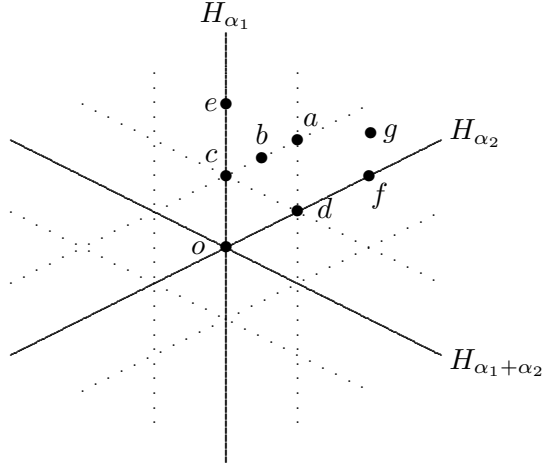
\centerline{{\bf Figure 4.1.}  Real parts of central characters in Table 4.1}

\medskip
\subsection{Tempered and square integrable representations.}

 The tempered (resp. square integrable) $\tilde H$-modules are the
ones which have the real parts of all their weights in the closure (resp.
interior) of the shaded region of Figure 4.2.  Let $t\in T$ be given by 
$t(X^{-\alpha_1}) = \pm q$, $t(X^{-\alpha_2})=\pm q$.  
This is a special case of the central character $t_b$ in Table 4.1. For this
particular special case there is one tempered representation with central
character $t$.
$$
\beginpicture
\setcoordinatesystem units <1cm,0.866cm>         
\setplotarea x from -4 to 4, y from -3 to 3.5    
\put{$H_{\alpha_1}$}[b] at 0 3.05
\put{$H_{\alpha_2}$}[l] at 2.7 1.6
\put{$H_{\alpha_1+\alpha_2}$}[l] at 2.7 -1.6
\put{$\bullet$} at -.8660 -1.5
\put{$s_2s_1s_2t_a$}[tr] at -.93 -1.4
\put{$\bullet$} at -.4330 -.75
\put{$t$}[tr] at -.46 -.85
\put{$\bullet$} at .4330 -.75
\put{$s_1t$}[l] at .6 -.8
\put{$\bullet$} at -.8660 0
\put{$s_2t$}[br] at -1 0.1
\put{$\bullet$} at 0 0
\put{$t_o$}[bl] at 0.1 0.25
\plot 0 3    0 -3 / 
\plot -2.5980 -1.5    2.5980 1.5 / 
\plot -2.5980  1.5    2.5980 -1.5 /  
\setdots
\plot  -.8660 -2.4340  -.8660 2.6340 /
\plot 2.5980 0.5    -1.7320 -2.0 /
\plot 1.7320  -2.0    -2.5980 0.5 /
\setdashes
\plot 0 0   -3 0 / 
\plot 0 0   1.5 -2.5980 /
\vshade -3 -2.8 0   0 -2.8 0   1.5 -2.8 -2.5980     /          
\endpicture
$$
\centerline{{\bf Figure 4.2.} Real parts of weights of tempered representations}

\medskip
The irreducible tempered representations with real central character are in 
one-to-one correspondence with the irreducible representations of the symmetric group $S_3$ (see
[BM]).  These representations are indexed by the partitions $(3), (21), (1^3)$
of $3$.   Equivalently, they can be indexed by the pairs
$(n,\rho)$ which appear in the Springer correspondence.  The $n$ and $\rho$ will
also be elements of the indexing triple for the corresponding tempered
representation of $\tilde H$.  Here $n$ is a
nilpotent element of the Lie algebra $\fg={\goth{sl}}_3$ and $\rho$ is an
irreducible representation of the component group $Z_G(n)/Z_G(n)^\circ$.
In type A the component group is always trivial.   For each root $\beta\in R$
let $e_\beta$ be an element of the root space $\fg_\beta$.  The three nilpotent
orbits in $\fg$ and the corresponding tempered representations of $\tilde H$ are
as in Table 4.2. 
$$\matrix{
\hbox{Nilpotent orbit} &Z_G(n)/Z_G(n)^\circ  &\hbox{Indexing triple}
&\hbox{Square integrable} &\hbox{$W$ representation}\cr
\hbox{regular} &1 &(t_a, e_{\alpha_1}+e_{\alpha_2},1) &\hbox{yes}
&(3)\cr
\hbox{subregular} &1  &(s_2s_1t,e_{\alpha_2},1)
&\hbox{no} &(21)\cr 
0  &1  &(t_o,0,1) &\hbox{no} &(1^3) \cr
}$$
\centerline{{\bf Table 4.2.}  Tempered representations and the Springer
correspondence}


$$
\beginpicture
\setcoordinatesystem units <1cm,0.866cm>         
\setplotarea x from -2 to 2, y from -1.5 to 1.5    
\put{$\bullet$} at .4330 .725
\put{$t_a$}[bl] at .53 .85
\put{$\bullet$} at -.4330 .725
\put{$s_1t_a$}[br] at -.53 .85
\put{$\bullet$} at -.8660 0
\put{$s_2s_1t_a$}[r] at -.96 0
\put{$\bullet$} at .8660 0
\put{$s_2t_a$}[l] at .96 0
\put{$\bullet$} at .4330 -.725
\put{$s_1s_2t_a$}[tl] at .53 -.85
\put{$\bullet$} at -.4330 -.725
\put{$s_2s_1s_2t_a$}[tr] at -.53 -.85
\plot -.4330 .725    -.8660 0 / 
\plot .8660 0  .4330 -.725 / 
\endpicture
\qquad
\beginpicture
\setcoordinatesystem units <1cm,0.866cm>         
\setplotarea x from -2 to 1.5, y from -1.5 to 1.5    
\put{$\bullet$} at .4330 .725
\put{$t_b$}[bl] at .53 .85
\put{$\bullet$} at -.4330 .725
\put{$s_1t_b$}[br] at -.53 .85
\put{$\bullet$} at -.8660 0
\put{$s_2s_1t_b$}[r] at -.96 0
\put{$\bullet$} at .8660 0
\put{$s_2t_b$}[l] at .96 0
\put{$\bullet$} at .4330 -.725
\put{$s_1s_2t_b$}[tl] at .53 -.85
\put{$\bullet$} at -.4330 -.725
\put{$s_2s_1s_2t_b$}[tr] at -.53 -.85
\plot .4330 .725   -.4330 .725 /  
\plot -.4330 .725    -.8660 0 / 
\plot .8660 0  .4330 -.725 / 
\plot .4330 -.725    -.4330 -.725 / 
\endpicture
\qquad
\beginpicture
\setcoordinatesystem units <1cm,0.866cm>         
\setplotarea x from 0 to 2, y from -1.5 to 1.5    
\put{$\bullet$} at .4330 .725
\put{$t_c$}[bl] at .53 .85
\put{$\bullet$} at .8660 0
\put{$s_2t_c$}[l] at .96 0
\put{$\bullet$} at .4330 -.725
\put{$s_1s_2t_c$}[tl] at .53 -.85
\endpicture
$$
$$
\beginpicture
\setcoordinatesystem units <1cm,0.866cm>         
\setplotarea x from -2 to 2, y from -1.5 to 1.5    
\put{$\bullet$} at .4330 .725
\put{$t_d$}[bl] at .53 .85
\put{$\bullet$} at -.4330 .725
\put{$s_1t_d$}[br] at -.53 .85
\put{$\bullet$} at -.8660 0
\put{$s_2s_1t_d$}[r] at -.96 0
\endpicture
\qquad
\beginpicture
\setcoordinatesystem units <1cm,0.866cm>         
\setplotarea x from -2 to 2, y from -1.5 to 1.5    
\put{$\bullet$} at .4330 .725
\put{$t_e$}[bl] at .53 .85
\put{$\bullet$} at .8660 0
\put{$s_2t_e$}[l] at .96 0
\put{$\bullet$} at .4330 -.725
\put{$s_1s_2t_e$}[tl] at .53 -.85
\plot  .4330 .725  .8660 0 / 
\plot .8660 0  .4330 -.725 / 
\endpicture
\qquad
\beginpicture
\setcoordinatesystem units <1cm,0.866cm>         
\setplotarea x from -2 to 1.5, y from -1.5 to 1.5    
\put{$\bullet$} at .4330 .725
\put{$t_f$}[bl] at .53 .85
\put{$\bullet$} at -.4330 .725
\put{$s_1t_f$}[br] at -.53 .85
\put{$\bullet$} at -.8660 0
\put{$s_2s_1t_f$}[r] at -.96 0
\plot .4330 .725   -.4330 .725 /  
\plot -.4330 .725    -.8660 0 / 
\endpicture
$$
$$
\beginpicture
\setcoordinatesystem units <1cm,0.866cm>         
\setplotarea x from -2 to 2, y from -1.5 to 1.5    
\put{$\bullet$} at .4330 .725
\put{$t_g$}[bl] at .53 .85
\put{$\bullet$} at -.4330 .725
\put{$s_1t_g$}[br] at -.53 .85
\put{$\bullet$} at -.8660 0
\put{$s_2s_1t_g$}[r] at -.96 0
\put{$\bullet$} at .8660 0
\put{$s_2t_g$}[l] at .96 0
\put{$\bullet$} at .4330 -.725
\put{$s_1s_2t_g$}[tl] at .53 -.85
\put{$\bullet$} at -.4330 -.725
\put{$s_2s_1s_2t_g$}[tr] at -.53 -.85
\plot .4330 .725   -.4330 .725 /  
\plot -.4330 .725    -.8660 0 / 
\plot -.8660 0    -.4330 -.725 / 
\plot  .4330 .725  .8660 0 / 
\plot .8660 0  .4330 -.725 / 
\plot .4330 -.725    -.4330 -.725 / 
\endpicture
\qquad
\beginpicture
\setcoordinatesystem units <1cm,0.866cm>         
\setplotarea x from -1 to 2, y from -1.5 to 1.5    
\put{$\bullet$} at .4330 .725
\put{$t_o$}[bl] at .53 .85
\endpicture
$$
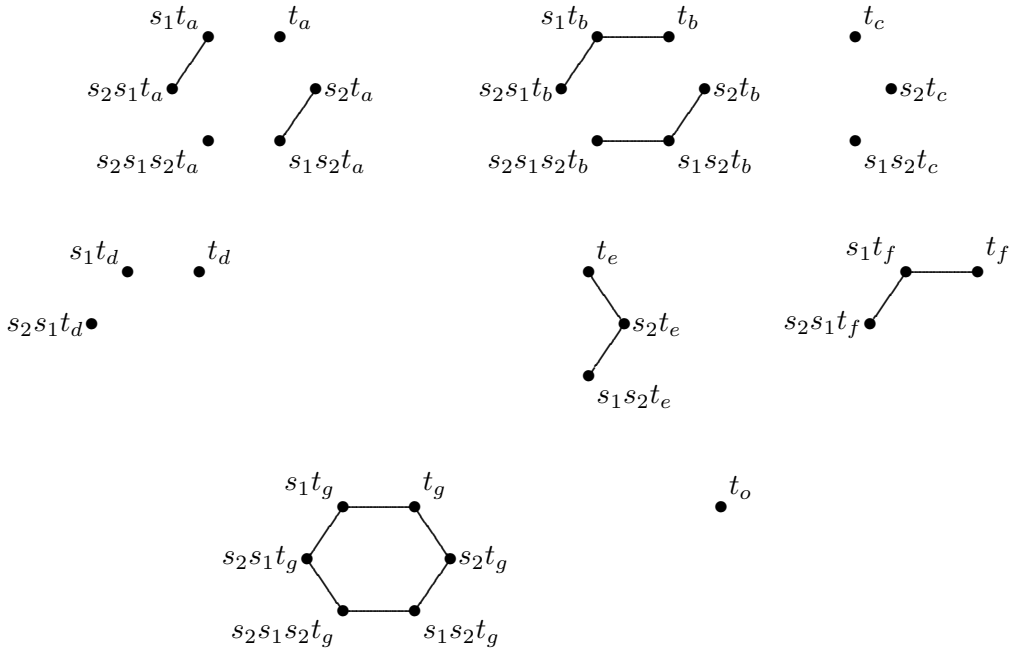
\centerline{{\bf Figure 4.3.}  Calibration graphs for central characters in Table
4.1}

\bigskip

\subsection{The analysis.}

\smallskip\noindent
{\it Central characters $t_a$, $t_b$, $t_g$:}  Since
$Z(t)=\emptyset$ these weights are regular.   Thus the representations
corresponding to these central characters are in one to one correspondence
with the connected components of the calibration graph $\Gamma(t)$ and
can be constructed explicitly with the use of Theorem 1.11.  Up to isomorphism
the principal series module $M(t_g)$ is the only irreducible $\tilde H$-module with
central character $t_g$.  The Langlands parameters for each module can be
determined from its weight structure and the indexing triple is then determined
from the Langlands data by using the induction theorem of Kazhdan and Lusztig (see
the discusssion in [BM, p.34]).

There is one special case of the central character $t_b$ when the irreducible
module constructed by applying Theorem 1.11 to the placed skew shape
$(t_b,\emptyset)$ is tempered.  This happens when $t_b =s_2s_1t$ for the weight
$t\in T$ given by $t(X^{-\alpha_1})=\pm q$, $t(X^{-\alpha_2})=\pm q$.
The indexing triple and the calibration set for this case are still given by
$(t_b,e_{\alpha_2},1)$ and $J=\{\alpha_2\}$, respectively.

\smallskip\noindent
{\it Central characters $t_c$ and $t_d$:}
One can use the defining relations of $\tilde H$ to check that the only
$1$-dimensional representations of $\tilde H$ are the ones with
central character $t_a$.   Construct two $3$-dimensional representations of 
$\tilde H$ by
$$\Ind_{\tilde H_{\{2\}}}^{\tilde H}(\CC v_c)
\qquad\hbox{and}\qquad
\Ind_{\tilde H_{\{2\}}}^{\tilde H}(\CC v_{s_2c}),$$
where $\CC v_c$ and $\CC v_{s_2c}$ are the two one-dimensional representations
of $\tilde H_{\{2\}}$ given by
$$
\matrix{
T_2v_c = qv_c, 
&&X^{\alpha_1}v_c = v_c ,
&&X^{\alpha_2}v_c = q^2v_c, \cr
\cr
T_2v_{s_2c} = -q^{-1}v_{s_2c} ,
&&X^{\alpha_1}v_{s_2c} = q^2v_{s_2c}, 
&&X^{\alpha_2}v_{s_2c} = q^{-2}v_{s_2c}. \cr
}$$
These representations must be irreducible since, if not, they would either have a
$1$-dimensional submodule or a one dimensional quotient.  But there are no
$1$-dimensional modules with central character $t_c$.

The central characters $t_c$ and $t_d$ are taken into each other
under the automorphism of the Dynkin diagram of $A_2$ which switches the
two nodes and thus these two central characters will produce modules
which have the same structure (up to twisting by the automorphism
which switches $\alpha_1$ and $\alpha_2$).  Thus the representations with
central character $t_d$ can be obtained from the ones with central character
$t_c$ by switching all $1$'s and $2$'s and changing all $c$'s to $d$'s.

\smallskip\noindent
{\it Central characters $t_e$ and $t_f$:}  Since $P(t_e)=\emptyset$, Kato's
irreducibility criterion (Theorem 1.16)  implies that the
principal series module $M(t_e)$ is irreducible.  By Theorem 1.15 this is the
only irreducible with central character $t_e$.  
As for the case of $t_c$ and $t_d$, the central characters 
$t_e$ and $t_f$ are taken into each other under the automorphism of the Dynkin
diagram of $A_2$.  Thus the irreducible representations with central character
$t_f$ can be obtained from the one with central character
$t_e$ by switching all $1$'s and $2$'s and changing all $e$'s to $f$'s.

\smallskip\noindent
{\it Central character $t_o$:}  Since $P(t_o)=\emptyset$, Kato's
irreducibility criterion (Theorem 1.16)  implies that the
principal series module $M(t_o)$ is irreducible. By Theorem 1.15 this is the
only irreducible with central character $t_o$.  

\section 5. Classification for $C_2$

\bigskip

The root system $R$ for $C_2$ has simple roots $\alpha_1$ and $\alpha_2$,
fundamental weights $\omega_1$ and $\omega_2$, and 
$$
\eqalign{
\langle \alpha_1,\alpha_2^\vee\rangle &= -2 \cr
\langle \alpha_2,\alpha_1^\vee\rangle &= -1, \cr}
\qquad\eqalign{
\omega_1 &= \alpha_1+\alpha_2\cr
\omega_2&=\hbox{$1\over2$}\alpha_1+\alpha_2, \cr }
\qquad\hbox{and}\qquad
\eqalign{
\alpha_1 &= 2\omega_1-2\omega_2 \cr 
\alpha_2&=-\omega_1+2\omega_2. \cr}
$$

\subsection{Irreducible representations.}

Table 5.1 lists the irreducible $\tilde H$-modules by their
central characters.  We have listed only those central characters $t$ for which
the principal series module $M(t)$ is not irreducible (see Theorem 1.16).  The
sets $P(t)$ and $Z(t)$ are as given in (1.7) and correspond to the choice of
representative for the central character displayed in Figure 5.1.  
 The Langlands parameters usually consist of a pair $(\cT,I)$
where $I$ is a subset of $\{1,2\}$ and $\cT$ is a
tempered representation for the parabolic subalgebra $\tilde H_I$.  In our cases
the tempered representation $\cT$ of $\tilde H_I$ is completely determined by a
character $t\in T$.  Specifically, $\cT$ is the only tempered representation of
$\tilde H_I$ which has $t$ as a weight. In the labeling in
Table 5.1 we have replaced the representation $\cT$ by the weight $t$.  The
notation for the nilpotent elements in the indexing triples is as in Table 5.2. 
For each calibrated module with central character
$t$ we have listed the subset $J\subseteq P(t)$ such that $(t,J)$ is the
corresponding placed skew shape (see Theorem 1.11).  The abbreviation
`nc' indicates modules that are not calibrated. 

$$
\matrix{
\hbox{Central} &P(t)  &\hbox{Dim.} &\hbox{Langlands} 
&\hbox{Indexing} &\hbox{Calibration} \cr
\hbox{char.}  &Z(t) & &\hbox{parameters} &\hbox{triple} &\hbox{set $J$}
\cr
\cr
t_a &\{\alpha_1,\alpha_2\} 
&1 &(s_1s_2s_1s_2t_a,\emptyset) &(t_a,0,1) &\emptyset\cr 
&\emptyset
&3 &(s_1t_a,\{1\}) &(t_a,e_{\alpha_1},1) &\{\alpha_1\}\cr
&&3 &(s_2t_a,\{2\}) &(t_a,e_{\alpha_2},1) &\{\alpha_2\}\cr
&&1 &\hbox{tempered} &(t_a,e_{\alpha_1}+e_{\alpha_2},1) 
&\{\alpha_1,\alpha_2\}\cr
\cr
t_b &\{\alpha_1,\alpha_1+\alpha_2, 
&3 &(t_b,\{2\})  &(t_b,0,1) &\hbox{nc}\cr 
&\qquad\alpha_1+2\alpha_2\} 
&1 &(s_1t_b,\{1\})  &(t_b,e_{\alpha_1},1)
&\{\alpha_1\}\cr  
&\{\alpha_2\}&1 &\hbox{tempered} &(t_b,e_{\alpha_1+\alpha_2},-1)
&\{\alpha_1,\alpha_1+\alpha_2\}\cr  
&&3 &\hbox{tempered} &(t_b,e_{\alpha_1+\alpha_2},1) 
&\hbox{nc}\cr
\cr
t_c &\{\alpha_1,\alpha_1+2\alpha_2\} 
&2 &(t_c,\{2\}) &(t_c,0,1) &\emptyset\cr  
&\emptyset &2 &(s_1t_c,\{1\})  &(t_c,e_{\alpha_1},1) &\{\alpha_1\}\cr 
&&2 &(s_1s_2t_c,\{1\})  &(t_c,e_{\alpha_1+2\alpha_2},1) 
&\{\alpha_1+2\alpha_2\}\cr 
&&2 &\hbox{tempered} &(t_c,e_{\alpha_1}+e_{\alpha_1+2\alpha_2},1) 
&\{\alpha_1,\alpha_1+2\alpha_2\}\cr
\cr
t_d &\{\alpha_2,\alpha_1+\alpha_2\} 
&4 &(t_d,\{1\}) &(t_d,0,1) &\hbox{nc}\cr 
&\{\alpha_1\} &4 &(s_2t_d,\{2\})&(t_d,e_{\alpha_2},1) &\hbox{nc}\cr
\cr
t_e &\{\alpha_1\} 
&4 &(s_2t_e,\{1\}) &(t_e,0,1) &\hbox{nc}\cr 
&\{\alpha_1+2\alpha_2\}&4 &\hbox{tempered}  &(t_e,e_{\alpha_1},1) &\hbox{nc}\cr
\cr
t_f &\{\alpha_1\} 
&4 &(t_f,\emptyset)  &(t_f,0,1) &\emptyset\cr 
&\emptyset &4 &(s_1t_f,\{1\}) &(t_f,e_{\alpha_1},1) &\{\alpha_1\}\cr
\cr
t_g &\{\alpha_2\} 
&4 &(t_g,\emptyset) 
&(t_g,0,1) &\emptyset\cr 
&\emptyset &4 &(s_2t_g,\{2\}) &(t_g,e_{\alpha_2},1) &\{\alpha_2\}\cr
\cr
}
$$
\centerline{{\bf Table 5.1.}  Irreducible (non principal series) representations}

\bigskip

Figure 5.1 displays the real parts of the central
characters in Table 5.1.  If $t\in T$ then the polar decomposition
$t=t_rt_c$ determines an element $\mu\in \RR^n$ such that
$t_r(X^\lambda) = e^{\langle \lambda, \mu \rangle}$ (see (1.3)).  
For each central character $t_p$ the point labeled by $p$ in Figure 5.1 is
the graph of the corresponding $\mu_p\in \RR^n$.  Assume (for pictorial
convenience) that $q$ is a positive real number and let
$$
H_\beta = \{ x\in \RR^n \ |\ \langle \beta,x\rangle = 0 \},
\qquad\hbox{and}\qquad
H_{\beta\pm \delta} = \{ x\in \RR^n \ |\ \langle \beta,x\rangle = \ln(q^{\pm2})
\},
$$
for each positive root $\beta$.  The dotted lines display the (affine) hyperplanes
$H_{\beta\pm\delta}$.

\vfill\eject

$$
\beginpicture
\setcoordinatesystem units <1cm,1cm>         
\setplotarea x from -4 to 4, y from -3.5 to 3.5    
\put{$H_{\alpha_1}$}[b] at 0 3.1
\put{$H_{\alpha_2}$}[l] at 3.1 3.1
\put{$H_{\alpha_1+\alpha_2}$}[r] at -3.1 3.1
\put{$H_{\alpha_1+2\alpha_2}$}[l] at 3.1 0
\put{$\bullet$} at 1 3
\put{$a$}[l] at 1.1 3
\put{$\bullet$} at 1 1
\put{$b,c$}[bl] at 1.3 1.05
\put{$\bullet$} at 0 2
\put{$d$}[r] at -0.1 2.1
\put{$\bullet$} at 1 0
\put{$e$}[t] at  1 -0.1 
\put{$\bullet$} at 1 2
\put{$f$}[l] at  1.1 2 
\put{$\bullet$} at 0.5 2.5
\put{$g$}[br] at  0.4 2.6 
\plot -3 -3   3 3 /
\plot  3 -3  -3 3 /
\plot  0  3   0 -3 /
\plot  3  0  -3  0 /
\setdots
\putrule from 1 3.5 to 1 -3.5
\putrule from -1 3.5 to -1 -3.5
\putrule from 3.5 1 to -3.5  1
\putrule from 3.5 -1 to -3.5 -1
\plot 3.5 1.5  -1.5 -3.5 /
\plot -3.5 -1.5  1.5 3.5 /
\plot -3.5 1.5  1.5 -3.5 /
\plot 3.5 -1.5  -1.5 3.5 /
\endpicture
$$
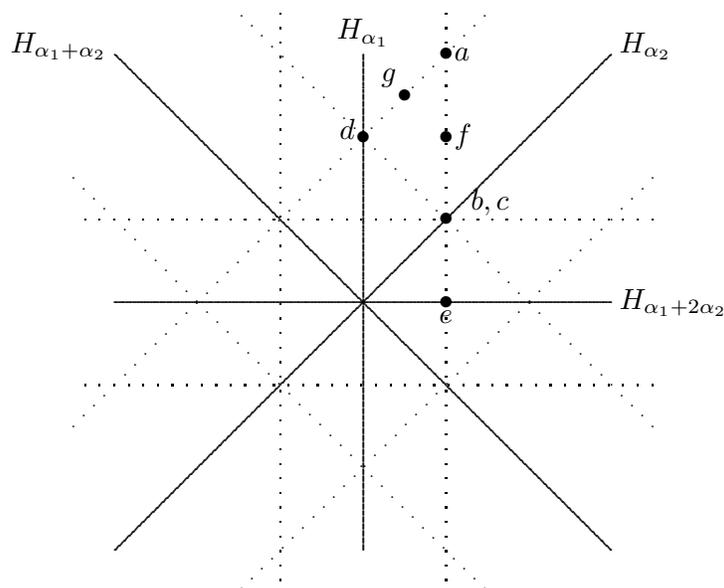
\centerline{{\bf Figure 5.1.}  Real parts of central characters in Table 5.1}

\subsection{Tempered and square integrable representations.}

The tempered (resp. square integrable) $\tilde H$-modules are the ones which have
the real parts of all their weights in the closure (resp. interior) of the
shaded region of Figure 5.2. 
$$
\beginpicture
\setcoordinatesystem units <1cm,1cm>         
\setplotarea x from -4 to 4, y from -4 to 4    
\put{$H_{\alpha_1}$}[b] at 0 3.1
\put{$H_{\alpha_2}$}[l] at 3.1 3.1
\put{$H_{\alpha_1+\alpha_2}$}[r] at -3.1 3.1
\put{$H_{\alpha_1+2\alpha_2}$}[l] at 3.1 0
\put{$\bullet$} at -1 -3
\put{$s_1s_2s_1s_2a$}[r] at -1.1 -3.1
\put{$\bullet$} at -1 -1
\put{$s_1s_2s_1b,s_1s_2s_1s_2c$}[tr] at -1.3 -.9
\put{$\bullet$} at 1 -1
\put{$s_1b$}[tl] at 1.2 -.9
\put{$\bullet$} at -1 0
\put{$s_1e$}[b] at  -1 0.1 
\put{$\bullet$} at 0 -1
\put{$s_2s_1e$} at  0.1 -1.2
\put{$\bullet$} at 0 0
\put{$t_o$}[bl] at  0.35 0.05
\plot -3 -3   3 3 /
\plot  3 -3  -3 3 /
\plot  0  3   0 -3 /
\plot  3  0  -3  0 /
\setdots
\putrule from 1 3.5 to 1 -3.5
\putrule from -1 3.5 to -1 -3.5
\putrule from 3.5 1 to -3.5  1
\putrule from 3.5 -1 to -3.5 -1
\plot 3.5 1.5  -1.5 -3.5 /
\plot -3.5 -1.5  1.5 3.5 /
\plot -3.5 1.5  1.5 -3.5 /
\plot 3.5 -1.5  -1.5 3.5 /
\vshade -3 -3 0   0 -3 0    3 -3 -3   /
\endpicture
$$
\centerline{{\bf Figure 5.2.}  Real parts of weights of tempered representations}

\vfill\eject

The irreducible tempered representations with real central character can be indexed
by the irreducible representations of the Weyl group $W$ of type $C_2$ (see
[BM, p. 34]).  Equivalently, these representations can be indexed by the pairs
$(n,\rho)$ which appear in the Springer correspondence.  The $n$ and $\rho$ will
also be elements of the indexing triple for the corresponding tempered
representation of $\tilde H$.   Here $n$ is a nilpotent 
element of the Lie algebra $\fg=\Lie(G)$, $G$ is the complex simple group
over $\CC$ of type $C_2$ and $\rho$ is an irreducible representation of the
component group $Z_G(n)/Z_G(n)^\circ$ (see [Ca]).  For each root $\beta\in R$ let
$e_\beta$ be an element of the root space $\fg_\beta$.  The four nilpotent orbits
in $\fg$ and the corresponding tempered representations of $\tilde H$ are as in
Table 5.2. We have used the notation of Carter [Ca, p.424] to label the
irreducible representations of the Weyl group.
$$\matrix{
\hbox{Nilpotent orbit} &Z_G(n)/Z_G(n)^\circ  &\hbox{Indexing triple}
&\hbox{Square integrable} &\hbox{$W$ representation}\cr
\hbox{regular} &1 &(t_a, e_{\alpha_1}+e_{\alpha_2},1) &\hbox{yes}
&(\emptyset,11)\cr
\hbox{subregular} &\ZZ/2\ZZ  &(t_b,e_{\alpha_1+\alpha_2},1)
&\hbox{yes} &(1,1)\cr 
&  &(t_b, e_{\alpha_1+\alpha_2},-1) &\hbox{yes} &(\emptyset,2) \cr
\hbox{minimal}  &1  &(t_e,e_{\alpha_1},1)  &\hbox{no} &(11,\emptyset) \cr
0  &1  &(t_o,0,1) &\hbox{no} &(2,\emptyset) \cr
}$$
\centerline{{\bf Table 5.2.}  Tempered representations and the Springer
correspondence}
\bigskip\noindent
The only other tempered representation is the square integrable
representation
$(t_c, e_{\alpha_1}+e_{\alpha_1+2\alpha_2},1)$.   
This representation does not have real central character.  It is 
the representation constructed in [Lu3, 4.14, 4.23]. (In Lusztig's notation
it is the star of the representation corresponding to the graph
$\cG'\oplus\cG''$.)

\vskip-.2in

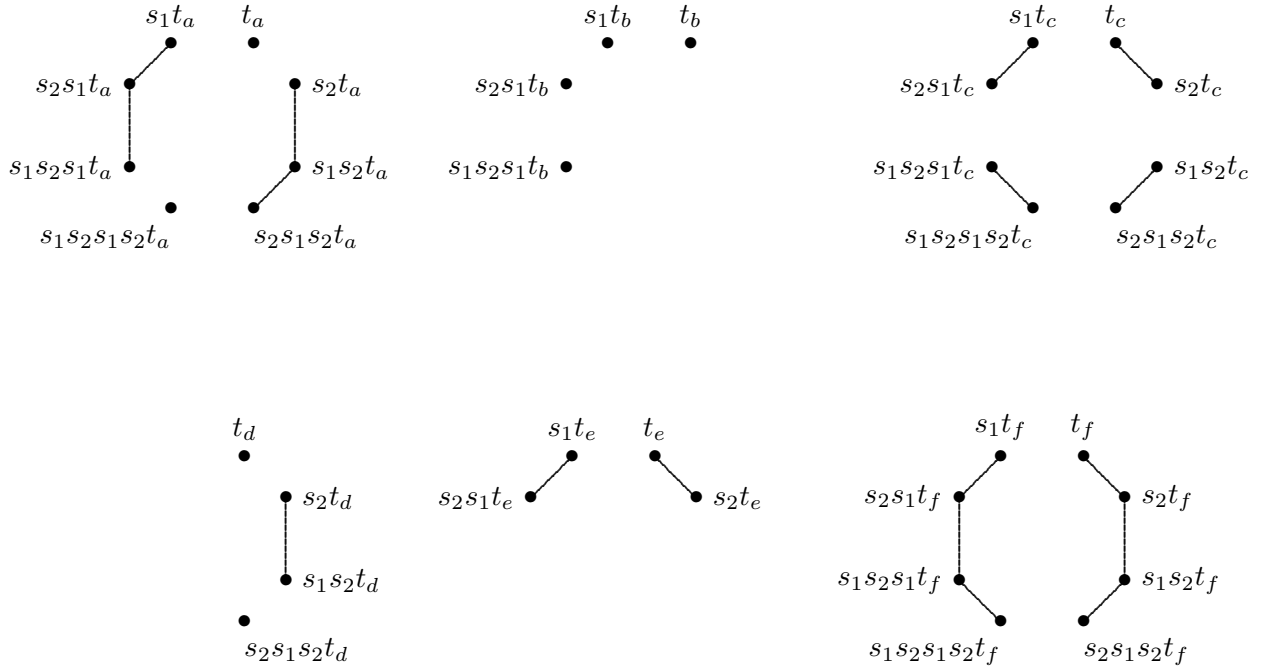
$$
\beginpicture
\setcoordinatesystem units <1cm,1cm>         
\setplotarea x from -2 to 2, y from -2 to 2    
\put{$\bullet$} at .5 1
\put{$t_a$}[b] at .5 1.2
\put{$\bullet$} at -.5 1
\put{$s_1t_a$}[b] at -.5 1.2
\put{$\bullet$} at .5 -1
\put{$s_2s_1s_2t_a$}[tl] at .5 -1.2
\put{$\bullet$} at -.5 -1
\put{$s_1s_2s_1s_2t_a$}[tr] at -.5 -1.2
\put{$\bullet$} at 1 .5
\put{$s_2t_a$}[l] at 1.2 .5
\put{$\bullet$} at -1 .5
\put{$s_2s_1t_a$}[r] at -1.2 .5
\put{$\bullet$} at 1 -.5
\put{$s_1s_2t_a$}[l] at  1.2 -.5 
\put{$\bullet$} at -1 -.5
\put{$s_1s_2s_1t_a$}[r] at  -1.2  -.5
\plot -.5 1   -1 .5 / 
\plot  -1 .5  -1 -.5 / 
\plot  1 .5  1 -.5 / 
\plot  1 -.5    .5 -1 / 
\endpicture
\qquad
\beginpicture
\setcoordinatesystem units <1cm,1cm>         
\setplotarea x from -2 to 2, y from -2 to 2    
\put{$\bullet$} at .5 1
\put{$t_b$}[b] at .5 1.2
\put{$\bullet$} at -.5 1
\put{$s_1t_b$}[b] at -.5 1.2
\put{$\bullet$} at -1 .5
\put{$s_2s_1t_b$}[r] at -1.2 .5
\put{$\bullet$} at -1 -.5
\put{$s_1s_2s_1t_b$}[r] at  -1.2  -.5
\endpicture
\qquad
\beginpicture
\setcoordinatesystem units <1cm,1cm>         
\setplotarea x from -2 to 2, y from -2 to 2    
\put{$\bullet$} at .5 1
\put{$t_c$}[b] at .5 1.2
\put{$\bullet$} at -.5 1
\put{$s_1t_c$}[b] at -.5 1.2
\put{$\bullet$} at -1 .5
\put{$s_2s_1t_c$}[r] at -1.2 .5
\put{$\bullet$} at -1 -.5
\put{$s_1s_2s_1t_c$}[r] at  -1.2  -.5
\put{$\bullet$} at 1 .5
\put{$s_2t_c$}[l] at 1.2 .5
\put{$\bullet$} at 1 -.5
\put{$s_1s_2t_c$}[l] at  1.2 -.5 
\put{$\bullet$} at .5 -1
\put{$s_2s_1s_2t_c$}[tl] at .5 -1.2
\put{$\bullet$} at -.5 -1
\put{$s_1s_2s_1s_2t_c$}[tr] at -.5 -1.2
\plot -.5 1   -1 .5 / 
\plot  -1 -.5  -.5  -1 / 
\plot  .5 1   1 .5 / 
\plot  1 -.5    .5 -1 / 
\endpicture
$$
$$\beginpicture
\setcoordinatesystem units <1cm,1cm>         
\setplotarea x from -2 to 2, y from -2 to 2    
\put{$\bullet$} at .5 1
\put{$t_d$}[b] at .5 1.2
\put{$\bullet$} at 1 .5
\put{$s_2t_d$}[l] at 1.2 .5
\put{$\bullet$} at 1 -.5
\put{$s_1s_2t_d$}[l] at  1.2 -.5 
\put{$\bullet$} at .5 -1
\put{$s_2s_1s_2t_d$}[tl] at .5 -1.2
\plot  1 .5  1 -.5 / 
\endpicture
\qquad
\beginpicture
\setcoordinatesystem units <1cm,1cm>         
\setplotarea x from -2 to 2, y from -2 to 2    
\put{$\bullet$} at .5 1
\put{$t_e$}[b] at .5 1.2
\put{$\bullet$} at -.5 1
\put{$s_1t_e$}[b] at -.5 1.2
\put{$\bullet$} at -1 .5
\put{$s_2s_1t_e$}[r] at -1.2 .5
\put{$\bullet$} at 1 .5
\put{$s_2t_e$}[l] at 1.2 .5
\plot -.5 1   -1 .5 / 
\plot  .5 1   1 .5 / 
\endpicture
\qquad
\beginpicture
\setcoordinatesystem units <1cm,1cm>         
\setplotarea x from -2 to 2, y from -2 to 2    
\put{$\bullet$} at .5 1
\put{$t_f$}[b] at .5 1.2
\put{$\bullet$} at -.5 1
\put{$s_1t_f$}[b] at -.5 1.2
\put{$\bullet$} at -1 .5
\put{$s_2s_1t_f$}[r] at -1.2 .5
\put{$\bullet$} at -1 -.5
\put{$s_1s_2s_1t_f$}[r] at  -1.2  -.5
\put{$\bullet$} at 1 .5
\put{$s_2t_f$}[l] at 1.2 .5
\put{$\bullet$} at 1 -.5
\put{$s_1s_2t_f$}[l] at  1.2 -.5 
\put{$\bullet$} at .5 -1
\put{$s_2s_1s_2t_f$}[tl] at .5 -1.2
\put{$\bullet$} at -.5 -1
\put{$s_1s_2s_1s_2t_f$}[tr] at -.5 -1.2
\plot -.5 1   -1 .5 / 
\plot  -1 .5  -1 -.5 / 
\plot  -1 -.5  -.5  -1 / 
\plot  .5 1   1 .5 / 
\plot  1 .5  1 -.5 / 
\plot  1 -.5    .5 -1 / 
\endpicture
$$
$$
\beginpicture
\setcoordinatesystem units <1cm,1cm>         
\setplotarea x from -2 to 2, y from -2 to 2    
\put{$\bullet$} at .5 1
\put{$t_g$}[b] at .5 1.2
\put{$\bullet$} at -.5 1
\put{$s_1t_g$}[b] at -.5 1.2
\put{$\bullet$} at -1 .5
\put{$s_2s_1t_g$}[r] at -1.2 .5
\put{$\bullet$} at -1 -.5
\put{$s_1s_2s_1t_g$}[r] at  -1.2  -.5
\put{$\bullet$} at 1 .5
\put{$s_2t_g$}[l] at 1.2 .5
\put{$\bullet$} at 1 -.5
\put{$s_1s_2t_g$}[l] at  1.2 -.5 
\put{$\bullet$} at .5 -1
\put{$s_2s_1s_2t_g$}[tl] at .5 -1.2
\put{$\bullet$} at -.5 -1
\put{$s_1s_2s_1s_2t_g$}[tr] at -.5 -1.2
\plot .5 1  -.5 1 / 
\plot -.5 1   -1 .5 / 
\plot  -1 .5  -1 -.5 / 
\plot  1 .5  1 -.5 / 
\plot  1 -.5    .5 -1 / 
\plot  .5 -1   -.5 -1 / 
\endpicture
$$
\nobreak
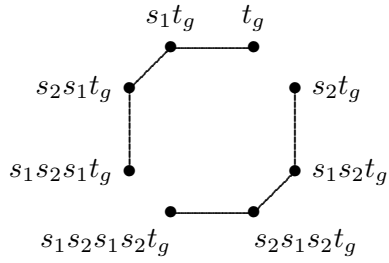
\centerline{{\bf Figure 5.3.}  Calibration graphs for central characters in
Table 5.1}

\bigskip

\subsection{The analysis.}

A general calculation with the defining relations of 
$\tilde H$ shows that there are only four one dimensional 
$\tilde H$-modules $\CC v_1$, $\CC v_2$, $\CC v_3$, and $\CC v_4$
which have weights $t_a$, $s_2s_1t_b$, $s_1t_b$ and $s_1s_2s_1s_2t_a$,
respectively.  These modules are given explicitly by
$$
\matrix{
T_1v_1=qv_1, &&T_2v_1=qv_1, 
&&X^{\alpha_1}v_1=q^2v_1, &&X^{\alpha_2}v_1=q^2v_1, \cr
T_1v_2=qv_2, &&T_2v_2=-q^{-1}v_2, 
&&X^{\alpha_1}v_2=q^2v_2,
&&X^{\alpha_2}v_2=q^{-2}v_2, \cr  
T_1v_3=-q^{-1}v_3, &&T_2v_3=qv_3,
&&X^{\alpha_1}v_3=q^{-2}v_3, &&X^{\alpha_2}v_3=q^2v_3,\cr 
T_1v_4=-q^{-1}v_4,  &&T_2v_4=-q^{-1}v_4, 
&&X^{\alpha_1}v_4=q^{-2}v_4,  &&X^{\alpha_2}v_4=q^{-2}v_4. \cr }$$

\smallskip\noindent
{\it Central character $t_a$:}  Since $Z(t_a)=\emptyset$, $t_a$ is regular
and thus all irreducible representations with central character $t_a$ are
calibrated.  They are in one to one correspondence with the connected components
of the calibration graph and can be constructed with the use of Theorem 1.11.
The Langlands parameters for each module can be determined from its
weight structure and the indexing triple
is then determined from the Langlands data by using the induction theorem of
Kazhdan-Lusztig (see the discussion in [BM, p.34]).

\smallskip\noindent
{\it Central character $t_b$:}  
We already know from our general computation above, that there are two
one-dimensional $\tilde H$-modules with central character $t_b$.
One has weight $s_1t_b$ and the other has weight $s_2s_1t_b$.
Let $\CC v_b$ and $\CC
v_{s_1b}$ be the one dimensional representations of $\tilde H_{\{1\}}$ given by
$$
\matrix{
T_1v_b=qv_b,  &&X^{\alpha_1}v_b=q^2v_b,
&&X^{\alpha_2}v_b=v_b, \cr 
T_1v_{s_1b}=-q^{-1}v_{s_1b}, &&X^{\alpha_1}v_{s_1b}=q^{-2} v_{s_1b},
&&X^{\alpha_2}v_{s_1b}=q^2v_{s_1b}. \cr }
$$
Let 
$$M_1 = \Ind_{\tilde H_{\{1\}}}^{\tilde H}(\CC v_b)
\qquad\hbox{and}\qquad
M_2 = \Ind_{\tilde H_{\{1\}}}^{\tilde H}(\CC v_{s_1b}).$$
By Lemma 1.17 these modules have weights
$\supp(M_1)=\{t_b,s_1t_b,s_2s_1t_b\}$ and
$\supp(M_2)=\{s_1t_b,$ $s_2s_1t_b, s_1s_2s_1t_b\}$, respectively.  Both
$M_1$ and $M_2$ are $4$ dimensional.
By Proposition 1.18 (c) one of the two operators $\tau_2\colon (M_1)_{s_2s_1t_b}\to
(M_1)_{s_1t_b}$ or $\tau_2\colon (M_1)_{s_1t_b}\to
(M_1)_{s_2s_1t_b}$ must have nonzero kernel.  This implies that
$M_1$ has either a $3$ dimensional submodule or a $3$ dimensional
quotient, call it $N_1$, with weights $\{t_b,s_1t_b\}$.  By Lemma 1.19,
any module $P$ with weights $\{t_b,s_1t_b\}$ must have
$\dim(P_{t_b}^{\rm gen})\ge 2$ and $\dim(P_{s_1t_b}^{\rm gen})\ge 1$.
It follows that $N_1$ is irreducible.  A similar argument
can be used to show that $M_2$ has either a $3$ dimensional submodule
or a $3$-dimensional quotient which must be irreducible.

The representation $N_1$ constructed in the
previous paragraph and the $1$ dimensional representation
with weight $s_1t_b$ are both tempered.  One obtains the
corresponding indexing triples  by comparing the Langlands
parameters for these modules with the labelings of the corresponding 
representations of $W$ in the Springer correspondence.  See [Ca, p.424],
[BM, p. 34] and Table 5.2.

\smallskip\noindent
{\it Central character $t_c$:}  Since $Z(t_c)=\emptyset$, $t_c$ is regular
and thus all irreducible representations with central character $t_c$ are
calibrated.  They are in one to one correspondence with the connected components
of the calibration graph and can be constructed with the use of Theorem 1.11.
The Langlands parameters for each module can be determined from its
weight structure.  The only representation for
which the indexing triple cannot be determined from the Langlands parameters
and the [KL] induction theorem (see the discussion in [BM, p. 34]) is the
tempered representation.  This representation is constructed in [Lu3, 4.14 and
4.23]. In Lusztig's notation, it is the star (see [Lu3,4.23]) of the representation
corresponding to the graph $\cG'\oplus\cG''$.
The indexing triple for this representation is given in the discussion for $B_2$ 
in [Lu3, 2.10].

\smallskip\noindent
{\it Central character $t_d$:}  
Let $\CC v_d$ and $\CC v_{s_2d}$ be the one dimensional
representations of $\tilde H_{\{2\}}$ given by
$$
\matrix{
T_2v_d=qv_d,  &&X^{\alpha_1}v_d=v_d,
&&X^{\alpha_2}v_d=q^2v_d, \cr 
T_2v_{s_2d}=-q^{-1}v_{s_2d}, &&X^{\alpha_1}v_{s_2d}=q^4 v_{s_2d},
&&X^{\alpha_2}v_{s_2d}=q^{-2}v_{s_2d}. \cr }
$$
Let 
$$M_1 = \Ind_{\tilde H_{\{2\}}}^{\tilde H}(\CC v_d)
\qquad\hbox{and}\qquad
M_2 = \Ind_{\tilde H_{\{2\}}}^{\tilde H}(\CC v_{s_2d}).$$
By Lemma 1.17 these modules have weights
$\supp(M_1)=\{t_d,s_2t_d,s_1s_2t_d\}$ and
$\supp(M_2)=\{s_2t_d,s_1s_2t_d,$ $s_2s_1s_2t_d\}$, respectively.  Both
$M_1$ and $M_2$ are $4$ dimensional.

Let $M$ be any $\tilde H$-module such that $M_{t_d}\ne 0$.
By Lemma 1.19 (a), $\dim(M_{t_d}^{\rm gen})\ge 2$.  Since
$\langle \alpha_2,\alpha_1^\vee\rangle \ne 0$ it follows from Lemma 1.19 (b) that
$\dim(M_{s_2t_d}^{\rm gen})\ge 1$.  Then, by Proposition 1.6, 
$\dim(M_{s_1s_2t_d}^{\rm gen})\ge 1$.  Adding these numbers up we see that
$\dim(M)\ge 4$.  It follows that $M_1$ is irreducible.
An analogous argument can be applied to conclude that $M_2$ is irreducible.

\smallskip\noindent
{\it Central character $t_e$:} 
Let $\CC v_e$ and $\CC v_{s_1e}$ be the one dimensional
representations of $\tilde H_{\{1\}}$ given by
$$
\matrix{
T_1v_e=qv_e,  &&X^{\alpha_1}v_e=q^2v_e,
&&X^{\alpha_2}v_e=q^{-1}v_e, \cr 
T_1v_{s_1e}=-q^{-1}v_{s_1e}, &&X^{\alpha_1}v_{s_1e}=q^{-2}v_{s_1e},
&&X^{\alpha_2}v_{s_1e}=qv_{s_1e}. \cr }
$$
Let 
$$M_1 = \Ind_{\tilde H_{\{1\}}}^{\tilde H}(\CC v_e)
\qquad\hbox{and}\qquad
M_2 = \Ind_{\tilde H_{\{1\}}}^{\tilde H}(\CC v_{s_1e}).$$
By Lemma 1.17 these modules have weights
$\supp(M_1)=\{t_e,s_2t_e\}$ and
$\supp(M_2)=\{s_1t_e,s_2s_1t_e\}$ respectively.  Both
$M_1$ and $M_2$ are $4$ dimensional.

Let $M$ be any $\tilde H$-module such that $M_{t_e}\ne 0$.
By Lemma 1.19 (a) and Proposition 1.6, $\dim(M_{t_e}^{\rm gen})
=\dim(M_{s_2t_e}^{\rm gen})\ge 2$. Thus $\dim(M)\ge 4$.  It follows that $M_1$ is
irreducible. An analogous argument can be applied to conclude that $M_2$ is
irreducible.

\smallskip\noindent
{\it Central character $t_f$ and $t_g$:}   These cases are handled 
in the same way as for the central character $t_a$.

\bigskip

\section 6. Classification for $G_2$

The root system $R$ for $G_2$ has simple roots $\alpha_1$ and $\alpha_2$,
fundamental weights $\omega_1$ and $\omega_2$, and 
$$
\eqalign{
\langle \alpha_1,\alpha_2^\vee\rangle &= -3 \cr
\langle \alpha_2,\alpha_1^\vee\rangle &= -1, \cr}
\qquad\eqalign{
\omega_1 &= 2\alpha_1+3\alpha_2\cr
\omega_2&=\alpha_1+2\alpha_2, \cr }
\qquad\hbox{and}\qquad
\eqalign{
\alpha_1 &= -\omega_1+2\omega_2 \cr 
\alpha_2&=2\omega_1-3\omega_2. \cr}
$$

\vfill\eject

\subsection{Irreducible representations.}

$$
\matrix{
\hbox{Central} &P(t)  &\hbox{Dimension} &\hbox{Langlands} 
&\hbox{Indexing} \cr
\hbox{character} &Z(t) &&\hbox{parameters} &\hbox{triple} \cr
\cr
t_a &\{\alpha_1,\alpha_2\}  
&1  &(t_a,\emptyset) &(t_a,0,1) \cr
&\emptyset &5  &(s_1t_a,\{1\}) &(t_a,e_{\alpha_1},1) \cr
&&5  &(s_2t_a,\{2\}) &(t_a,e_{\alpha_2},1) \cr
&&1  &\hbox{tempered} &(t_a,e_{\alpha_1}+e_{\alpha_2},1) \cr
\cr
t_b  &\{\alpha_1\} 
&6  &(t_b,\emptyset) &(t_b,0,1) \cr
&\emptyset &6 &(s_1t_b,\{1\}) &(t_b,e_{\alpha_1},1) \cr
\cr
t_c  &\{\alpha_1,\alpha_1+3\alpha_2\} 
&2  &(t_c,\{2\}) &(t_c,0,1) \cr  
&\emptyset &4 &(s_1t_c,\{1\}) &(t_c,e_{\alpha_1},1) \cr
&&4 &(s_1s_2t_c,\{1\}) &(t_c,e_{\alpha_1+3\alpha_2},1) \cr 
&&2 &\hbox{tempered}
&(t_c,e_{\alpha_1}+e_{\alpha_1+3\alpha_2},1)\cr
\cr
t_d &\{\alpha_1,\alpha_1+2\alpha_2\} 
&3 &(t_d,\{2\}) &(t_d,0,1)  \cr
&\emptyset &3 &(s_1t_d,\{1\})
&(t_d,e_{\alpha_1},1) \cr
&&3 &(s_2s_1s_2t_d,\{2\}) &(t_d,e_{\alpha_1+2\alpha_2},1)\cr 
&&3 &\hbox{tempered} &(t_d,e_{\alpha_1}+e_{\alpha_1+2\alpha_2},1) \cr
\cr
t_e &\{\alpha_1,\alpha_1+2\alpha_2,\hfill 
&3 &(t_e,\{2\}) &(t_e,0,1) \cr
&\alpha_1+\alpha_2,\alpha_1+3\alpha_2\} 
&1 &(s_1t_e,\{1\})  &(t_e,e_{\alpha_1},1) \cr 
&\{\alpha_2\} &2 &(s_2s_1t_e,\{2\}) &(t_e,e_{\alpha_1+\alpha_2},1) \cr
&&1  &\hbox{tempered} &(t_e,e_{\alpha_1}+e_{\alpha_1+2\alpha_2},(21)) \cr  
&&3  &\hbox{tempered} &(t_e,e_{\alpha_1}+e_{\alpha_1+2\alpha_2},(3)) \cr
\cr
t_f &\{\alpha_1,2\alpha_1+3\alpha_2\}  
&6 &(t_f,\{1\}) &(t_f,0,1) \cr
&\{\alpha_1+3\alpha_2\} &6 &(s_1t_f,\{1\}) &(t_f,e_{\alpha_1},1) \cr
\cr
t_g &\{\alpha_1\}  
&6 &(t_g,\{2\}) &(t_g,0,1) \cr
&\{\alpha_1+2\alpha_2\} &6 &\hbox{tempered} &(t_g,e_{\alpha_1},1)\cr
\cr
t_h &\{\alpha_2\} 
&6 &(t_h,\emptyset) &(t_h,0,1) \cr
&\emptyset &6 &(s_2t_h,\{2\}) &(t_h,e_{\alpha_2},1) \cr
\cr
t_i &\{\alpha_2,\alpha_1+\alpha_2\}  
&6 &(t_i,\{1\}) &(t_i,0,1) \cr
&\{\alpha_1\} &6 &(s_2t_i,\{2\}) &(t_i,e_{\alpha_2},1) \cr
\cr
t_j &\{\alpha_2\}  
&6 &(t_j,\{1\}) &(t_j,0,1) \cr
&\{2\alpha_1+3\alpha_2\} &6 &\hbox{tempered} &(t_j,e_{\alpha_2},1) \cr
}
$$
\centerline{{\bf Table 6.1.}  Irreducible (non principal series) representations}

\bigskip
\noindent
Table 6.1 lists the irreducible
$\tilde H$-modules by their central characters.  We have listed only those central
characters $t$ for which the principal series module $M(t)$ is not irreducible
(see Theorem 1.16).  The sets $P(t)$ and $Z(t)$ are as given in (1.7) and
correspond to the choice of representative for the central character displayed in
Figure 6.1.   
 The Langlands parameters usually consist of a pair $(\cT,I)$
where $I$ is a subset of $\{1,2\}$ and $\cT$ is a
tempered representation for the parabolic subalgebra $\tilde H_I$.  In our cases
the tempered representation $\cT$ of $\tilde H_I$ is completely determined by a
character $t\in T$.  Specifically, $\cT$ is the only tempered representation of
$\tilde H_I$ which has $t$ as a weight. In the labeling in
Table 6.1 we have replaced the representation $\cT$ by the weight $t$.  The
notation for the nilpotent elements in the indexing triples is as in Table 6.3. 

Table 6.2 lists the irreducible calibrated $\tilde H$-modules.  For each module
with central character $t$ we have listed the subset $J\subseteq P(t)$ such that
$(t,J)$ is the corresponding placed skew shape (see Theorem 1.11).  We have listed
only those central  characters $t$ for which the principal series module $M(t)$ is
not irreducible (see Theorem 1.16).
$$
\matrix{
\hbox{Central} &P(t) &Z(t) &\hbox{Dimension}  
&\hbox{Indexing} &\hbox{Calibration} \cr
\hbox{character} &&& &\hbox{triple} &\hbox{set $J$} \cr
\cr
t_a &\{\alpha_1,\alpha_2\} &\emptyset 
&1  &(t_a,0,1) &\emptyset\cr
&&&5  &(t_a,e_{\alpha_1},1) &\{\alpha_1\}\cr
&&&5  &(t_a,e_{\alpha_2},1) &\{\alpha_2\}\cr
&&&1  &(t_a,e_{\alpha_1}+e_{\alpha_2},1) &\{\alpha_1,\alpha_2\}\cr
\cr
t_b  &\{\alpha_1\} &\emptyset
&6  &(t_b,0,1) &\emptyset\cr
&&&6  &(t_b,e_{\alpha_1},1) &\{\alpha_1\}\cr
\cr
t_c  &\{\alpha_1,\alpha_1+3\alpha_2\} &\emptyset
&2  &(t_c,0,1) &\emptyset\cr
&&&4  &(t_c,e_{\alpha_1},1) &\{\alpha_1\}\cr
&&&4  &(t_c,e_{\alpha_1+3\alpha_2},1) &\{\alpha_1+3\alpha_2\}\cr
&&&2  &(t_c,e_{\alpha_1}+e_{\alpha_1+3\alpha_2},1) 
&\{\alpha_1,\alpha_1+3\alpha_2\}\cr
\cr
t_d &\{\alpha_1,\alpha_1+2\alpha_2\} &\emptyset
&3 &(t_d,0,1) &\emptyset \cr
&&&3  &(t_d,e_{\alpha_1},1) &\{\alpha_1\}\cr
&&&3  &(t_d,e_{\alpha_1+2\alpha_2},1) &\{\alpha_1+2\alpha_2\}\cr 
&&&3  &(t_d,e_{\alpha_1}+e_{\alpha_1+2\alpha_2},1) 
&\{\alpha_1,\alpha_1+2\alpha_2\}\cr
\cr
t_e &\{\alpha_1,\alpha_1+2\alpha_2,\hfill &\{\alpha_2\}
&1  &(t_e,e_{\alpha_1},1)
&\{\alpha_1\}\cr 
&\alpha_1+\alpha_2,\alpha_1+3\alpha_2\} 
&&2  &(t_e,e_{\alpha_1+\alpha_2},1)
&\{\alpha_1,\alpha_1+\alpha_2\}\cr 
&&&1  &(t_e,e_{\alpha_1}+e_{\alpha_1+2\alpha_2},(21)) 
&P(t_e)\setminus \{\alpha_1+3\alpha_2\}\cr 
\cr
t_h &\{\alpha_2\}  &\emptyset
&6 &(t_h,0,1) &\emptyset \cr
&&&6 &(t_h,e_{\alpha_2},1) &\{\alpha_2\} \cr
}
$$
\bigskip
\centerline{{\bf Table 6.2.}  Calibrated irreducible (non principal series)
representations}

\bigskip
\bigskip

Figure 6.1 displays the real parts of the central
characters in Table 6.1.  If $t\in T$ then the polar decomposition
$t=t_rt_c$ determines an element $\nu\in \RR^n$ such that
$t_r(X^\lambda) = e^{\langle \nu,\lambda\rangle}$ (see (1.3)).  
For each central character $t_p$ the point labeled by $p$ in Figure 6.1 is
the graph of the corresponding $\nu_p\in \RR^n$.  Assume (for pictorial
convenience) that $q$ is a positive real number and let
$$
H_\beta = \{ x\in \RR^n \ |\ \langle \beta,x\rangle = 0 \},
\qquad\hbox{and}\qquad
H_{\beta\pm \delta} = \{ x\in \RR^n \ |\ \langle \beta,x\rangle = \ln(q^{\pm2})
\},
$$
for each positive root $\beta$.  The dotted lines display the (affine) hyperplanes
$H_{\beta\pm\delta}$.
\vfill\eject

$$
\beginpicture
\setcoordinatesystem units <1cm,1cm>         
\setplotarea x from -5 to 6, y from -5 to 5    
\put{$H_{\alpha_1}$}[b] at 0 5.1
\put{$H_{\alpha_2}$}[b] at 2.309 4.1
\put{$H_{\alpha_1+\alpha_2}$}[b] at -2.309 4.1
\put{$H_{\alpha_1+2\alpha_2}$}[tl] at 4.5 -0.1
\put{$H_{\alpha_1+3\alpha_2}$}[tl] at 3.8 2.8
\put{$H_{2\alpha_1+3\alpha_2}$}[tr] at -3.5 2.8
\put{$\bullet$} at .75 3.897
\put{$a$}[l] at .85 3.897
\put{$\bullet$} at .75 2.598
\put{$b$}[l] at .85 2.598
\put{$\bullet$} at .75 1.299
\put{$c,d,e$}[l] at .85 1.25
\put{$\bullet$} at .75 .4330
\put{$f$}[tl] at .85 .44
\put{$\bullet$} at .75 0
\put{$g$}[tl] at .80  -0.05 
\put{$\bullet$} at .375 3.2475
\put{$h$}[r] at .375 3.5
\put{$\bullet$} at 0 2.598
\put{$i$}[r] at -0.15 2.598
\put{$\bullet$} at -1.125 0.6495
\put{$j$}[b] at -1.15 0.75
\putrule from 0 5 to 0 -5        
\plot -2.309 -4   2.309 4 /      
\putrule from 5 0  to -5 0       
\plot -4 -2.309  4 2.309 /   
\plot -4 2.309  4 -2.309 /   
\plot -2.309 4  2.309 -4 /   
\setdots
\plot .75 4.9   .75 -4.9 /  
\plot -.75 4.9   -.75 -4.9 /  
\plot -1.25 4.763  3.5 -3.464 /   
\plot 1.25 4.763  -3.5 -3.464 /   
\plot -1.25 -4.763  3.5 3.464 /   
\plot 1.25 -4.763  -3.5 3.464 /   
\plot  4 -1.443 -4 3.175 /   
\plot  -4 1.443 4 -3.175 /   
\plot  -4 -1.443  4 3.175 /   
\plot  4 1.443  -4 -3.175 /   
\plot  -5 1.299  5 1.299 /   
\plot  -5 -1.299  5 -1.299 /   
\endpicture
$$
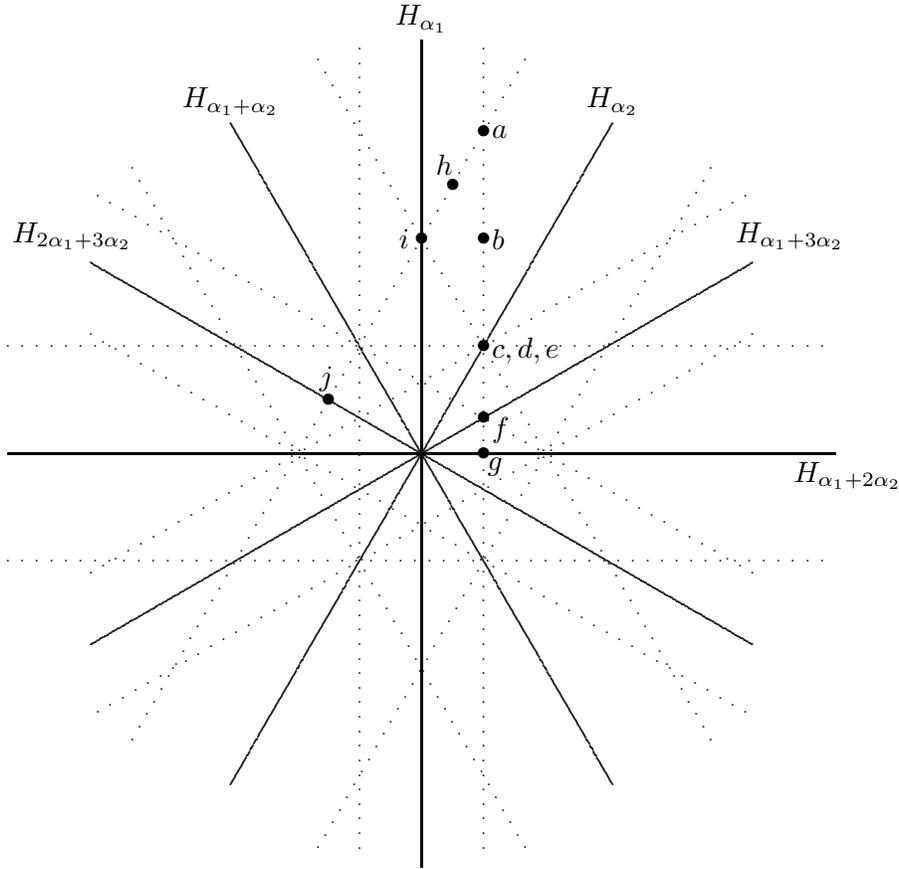
\centerline{{\bf Figure 6.1.}  Real parts of central characters in Table 6.1}

\subsection{Tempered and square integrable representations.}

The irreducible tempered representations with real central character can be indexed
by the irreducible representations of the Weyl group $W$ of type $G_2$ (see
[BM, p. 34]).  Equivalently, these representations can be indexed by the pairs
$(n,\rho)$ which appear in the Springer correspondence.  The $n$ and $\rho$ will
also be elements of the indexing triple for the corresponding tempered
representation of $\tilde H$.   Here $n$ is a nilpotent 
element of the Lie algebra $\fg=\Lie(G)$, $G$ is the complex simple group
over $\CC$ of type $G_2$ and $\rho$ is an irreducible representation of the
component group $Z_G(n)/Z_G(n)^\circ$ (see [Ca]).  For each root $\beta\in R$ let
$e_\beta$ be an element of the root space $\fg_\beta$. The five nilpotent orbits in
$\fg$ and the corresponding tempered representations of $\tilde H$ are as in Table
6.3. The notation $S_3$ denotes the symmetric group on three elements, which has
irreducible representations indexed by the partitions $(3), (21), (1^3)$ of $3$.
We have used the notation of Carter [Ca, p.427] to label the irreducible
representations of the Weyl group $W$ of type $G_2$.
$$\matrix{
\hbox{Nilpotent orbit} &&Z_G(n)/Z_G(n)^o  
&&\hbox{Indexing
triple} &&\hbox{Sq. int.} &&\hbox{$W$ rep.}\cr
\hbox{regular} &&1 &&(t_a, e_{\alpha_1}+e_{\alpha_2},1) &&\hbox{yes}
&&\phi_{1,0}\cr
\hbox{subregular} &&S_3  &&(t_e, e_{\alpha_1}+e_{\alpha_1+2\alpha_2},(3))
&&\hbox{yes} &&\phi_{2,1} \cr 
&&  &&(t_e, e_{\alpha_1}+e_{\alpha_1+2\alpha_2},(21)) &&\hbox{yes} 
&&{\phi_{1,3}}'\cr
\hbox{subminimal} &&1 &&(t_j,e_{\alpha_2},1) &&\hbox{no} 
&&\phi_{2,2} \cr
\hbox{minimal}  &&1  &&(t_g,e_{\alpha_1},1)  &&\hbox{no} 
&&{\phi_{1,3}}'' \cr
0  &&1  &&(t_o,0,1) &&\hbox{no}  &&\phi_{1,6} \cr
}$$
\centerline{{\bf Table 6.3.}  Tempered representations and the Springer
correspondence}

\vfill\eject

\medskip\noindent
The only other tempered representations are the representations
labeled by the triples $(t_c, e_{\alpha_1}+e_{\alpha_1+3\alpha_2},1)$ and
$(t_d, e_{\alpha_1}+e_{\alpha_1+2\alpha_2},1)$.  
These representations are square integrable but do not have real central character.

The modules labeled by $(t_c,e_{\alpha_1}+e_{\alpha_1+3\alpha_2},1)$,
$(t_d,e_{\alpha_1}+e_{\alpha_1+2\alpha_2},1)$,
$(t_e,e_{\alpha_1}+e_{\alpha_1+2\alpha_2},(3))$,
$(t_e,e_{\alpha_1}+e_{\alpha_1+2\alpha_2},(21))$
are the ones constructed by Lusztig in [Lu3] 4.20, 4.19, 4.7 and 4.22
respectively.  In Lusztig's notation these are the stars (see [Lu3, 4.23]) of
the modules labeled by the graphs $\cG''$, $\cG'$, $\cG$ and $\cG'''$,
respectively.

$$
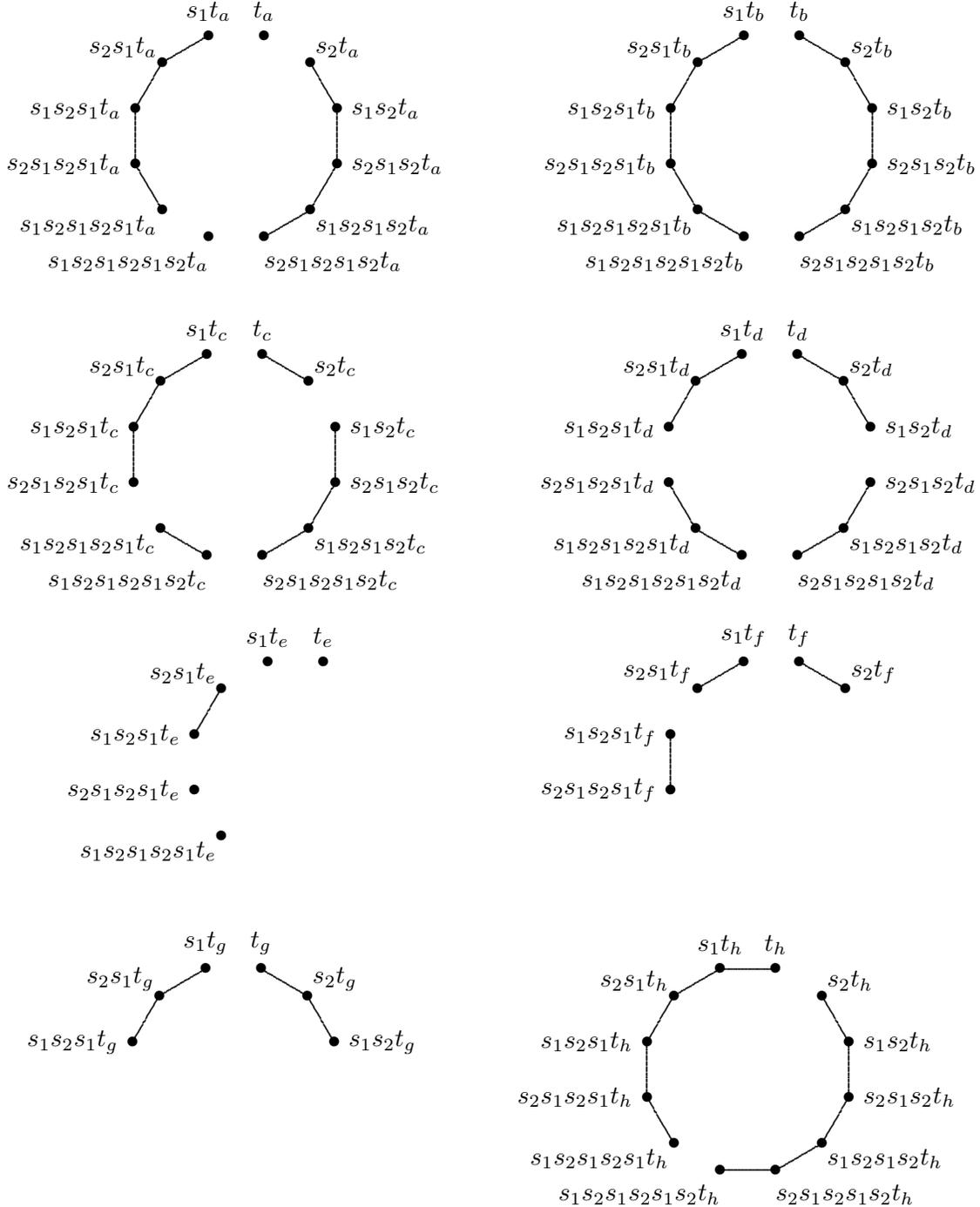
\beginpicture
\setcoordinatesystem units <1cm,1cm>         
\setplotarea x from -2 to 2, y from -2 to 2    
\put{$\bullet$} at .385 1.399
\put{$t_a$}[b] at .385 1.6
\put{$\bullet$} at -.385 1.399
\put{$s_1t_a$}[b] at -.385 1.6
\put{$\bullet$} at -1.026 1.026
\put{$s_2s_1t_a$}[br] at -1.1 1.1
\put{$\bullet$} at  -1.399 .385
\put{$s_1s_2s_1t_a$}[r] at  -1.6 .385
\put{$\bullet$} at -1.399 -.385
\put{$s_2s_1s_2s_1t_a$}[r] at -1.6 -.385
\put{$\bullet$} at -1.026 -1.026
\put{$s_1s_2s_1s_2s_1t_a$}[tr] at -1.1 -1.1
\put{$\bullet$} at 1.026 1.026
\put{$s_2t_a$}[bl] at 1.1 1.1
\put{$\bullet$} at  1.399 .385
\put{$s_1s_2t_a$}[l] at  1.6 .385
\put{$\bullet$} at  1.399 -.385
\put{$s_2s_1s_2t_a$}[l] at  1.6 -.385
\put{$\bullet$} at 1.026 -1.026
\put{$s_1s_2s_1s_2t_a$}[tl] at 1.1 -1.1
\put{$\bullet$} at .385 -1.399
\put{$s_2s_1s_2s_1s_2t_a$}[tl] at .385 -1.6
\put{$\bullet$} at -.385 -1.399
\put{$s_1s_2s_1s_2s_1s_2t_a$}[tr] at -.385 -1.6
\plot -1.026 1.026   -.385  1.399 / 
\plot -1.026 1.026   -1.399 .385 / 
\plot -1.399 -.385   -1.399 .385 / 
\plot -1.399 -.385   -1.026 -1.026 / 
\plot 1.026 1.026  1.399 .385 / 
\plot 1.399 -.385   1.399 .385  / 
\plot 1.399 -.385   1.026 -1.026  / 
\plot .385 -1.399   1.026 -1.026  / 
\endpicture
\qquad\qquad
\beginpicture
\setcoordinatesystem units <1cm,1cm>         
\setplotarea x from -2 to 2, y from -2 to 2    
\put{$\bullet$} at .385 1.399
\put{$t_b$}[b] at .385 1.6
\put{$\bullet$} at -.385 1.399
\put{$s_1t_b$}[b] at -.385 1.6
\put{$\bullet$} at -1.026 1.026
\put{$s_2s_1t_b$}[br] at -1.1 1.1
\put{$\bullet$} at  -1.399 .385
\put{$s_1s_2s_1t_b$}[r] at  -1.6 .385
\put{$\bullet$} at -1.399 -.385
\put{$s_2s_1s_2s_1t_b$}[r] at -1.6 -.385
\put{$\bullet$} at -1.026 -1.026
\put{$s_1s_2s_1s_2s_1t_b$}[tr] at -1.1 -1.1
\put{$\bullet$} at 1.026 1.026
\put{$s_2t_b$}[bl] at 1.1 1.1
\put{$\bullet$} at  1.399 .385
\put{$s_1s_2t_b$}[l] at  1.6 .385
\put{$\bullet$} at  1.399 -.385
\put{$s_2s_1s_2t_b$}[l] at  1.6 -.385
\put{$\bullet$} at 1.026 -1.026
\put{$s_1s_2s_1s_2t_b$}[tl] at 1.1 -1.1
\put{$\bullet$} at .385 -1.399
\put{$s_2s_1s_2s_1s_2t_b$}[tl] at .385 -1.6
\put{$\bullet$} at -.385 -1.399
\put{$s_1s_2s_1s_2s_1s_2t_b$}[tr] at -.385 -1.6
\plot -1.026 1.026   -.385  1.399 / 
\plot -1.026 1.026   -1.399 .385 / 
\plot -1.399 -.385   -1.399 .385 / 
\plot -1.399 -.385   -1.026 -1.026 / 
\plot -.385 -1.399    -1.026 -1.026 / 
\plot .385 1.399   1.026 1.026 / 
\plot 1.026 1.026  1.399 .385 / 
\plot 1.399 -.385   1.399 .385  / 
\plot 1.399 -.385   1.026 -1.026  / 
\plot .385 -1.399   1.026 -1.026  / 
\endpicture
$$
$$
\beginpicture
\setcoordinatesystem units <1cm,1cm>         
\setplotarea x from -2 to 2, y from -2 to 2    
\put{$\bullet$} at .385 1.399
\put{$t_c$}[b] at .385 1.6
\put{$\bullet$} at -.385 1.399
\put{$s_1t_c$}[b] at -.385 1.6
\put{$\bullet$} at -1.026 1.026
\put{$s_2s_1t_c$}[br] at -1.1 1.1
\put{$\bullet$} at  -1.399 .385
\put{$s_1s_2s_1t_c$}[r] at  -1.6 .385
\put{$\bullet$} at -1.399 -.385
\put{$s_2s_1s_2s_1t_c$}[r] at -1.6 -.385
\put{$\bullet$} at -1.026 -1.026
\put{$s_1s_2s_1s_2s_1t_c$}[tr] at -1.1 -1.1
\put{$\bullet$} at 1.026 1.026
\put{$s_2t_c$}[bl] at 1.1 1.1
\put{$\bullet$} at  1.399 .385
\put{$s_1s_2t_c$}[l] at  1.6 .385
\put{$\bullet$} at  1.399 -.385
\put{$s_2s_1s_2t_c$}[l] at  1.6 -.385
\put{$\bullet$} at 1.026 -1.026
\put{$s_1s_2s_1s_2t_c$}[tl] at 1.1 -1.1
\put{$\bullet$} at .385 -1.399
\put{$s_2s_1s_2s_1s_2t_c$}[tl] at .385 -1.6
\put{$\bullet$} at -.385 -1.399
\put{$s_1s_2s_1s_2s_1s_2t_c$}[tr] at -.385 -1.6
\plot -1.026 1.026   -.385  1.399 / 
\plot -1.026 1.026   -1.399 .385 / 
\plot -1.399 -.385   -1.399 .385 / 
\plot -1.026 -1.026   -.385 -1.399  / 
\plot .385 1.399   1.026 1.026 / 
\plot 1.399 -.385   1.399 .385  / 
\plot 1.399 -.385   1.026 -1.026  / 
\plot .385 -1.399   1.026 -1.026  / 
\endpicture
\qquad\qquad
\beginpicture
\setcoordinatesystem units <1cm,1cm>         
\setplotarea x from -2 to 2, y from -2 to 2    
\put{$\bullet$} at .385 1.399
\put{$t_d$}[b] at .385 1.6
\put{$\bullet$} at -.385 1.399
\put{$s_1t_d$}[b] at -.385 1.6
\put{$\bullet$} at -1.026 1.026
\put{$s_2s_1t_d$}[br] at -1.1 1.1
\put{$\bullet$} at  -1.399 .385
\put{$s_1s_2s_1t_d$}[r] at  -1.6 .385
\put{$\bullet$} at -1.399 -.385
\put{$s_2s_1s_2s_1t_d$}[r] at -1.6 -.385
\put{$\bullet$} at -1.026 -1.026
\put{$s_1s_2s_1s_2s_1t_d$}[tr] at -1.1 -1.1
\put{$\bullet$} at 1.026 1.026
\put{$s_2t_d$}[bl] at 1.1 1.1
\put{$\bullet$} at  1.399 .385
\put{$s_1s_2t_d$}[l] at  1.6 .385
\put{$\bullet$} at  1.399 -.385
\put{$s_2s_1s_2t_d$}[l] at  1.6 -.385
\put{$\bullet$} at 1.026 -1.026
\put{$s_1s_2s_1s_2t_d$}[tl] at 1.1 -1.1
\put{$\bullet$} at .385 -1.399
\put{$s_2s_1s_2s_1s_2t_d$}[tl] at .385 -1.6
\put{$\bullet$} at -.385 -1.399
\put{$s_1s_2s_1s_2s_1s_2t_d$}[tr] at -.385 -1.6
\plot -1.026 1.026   -.385  1.399 / 
\plot -1.026 1.026   -1.399 .385 / 
\plot -1.399 -.385   -1.026 -1.026 / 
\plot -.385 -1.399  -1.026 -1.026 / 
\plot .385 1.399   1.026 1.026 / 
\plot 1.026 1.026  1.399 .385 / 
\plot 1.399 -.385   1.026 -1.026  / 
\plot .385 -1.399   1.026 -1.026  / 
\endpicture
$$
$$
\beginpicture
\setcoordinatesystem units <1cm,1cm>         
\setplotarea x from -2 to 2, y from -2 to 2    
\put{$\bullet$} at .385 1.399
\put{$t_e$}[b] at .385 1.6
\put{$\bullet$} at -.385 1.399
\put{$s_1t_e$}[b] at -.385 1.6
\put{$\bullet$} at  -1.399 .385
\put{$s_1s_2s_1t_e$}[r] at  -1.6 .385
\put{$\bullet$} at -1.399 -.385
\put{$s_2s_1s_2s_1t_e$}[r] at -1.6 -.385
\put{$\bullet$} at -1.026 1.026
\put{$s_2s_1t_e$}[br] at -1.1 1.1
\put{$\bullet$} at -1.026 -1.026
\put{$s_1s_2s_1s_2s_1t_e$}[tr] at -1.1 -1.1
\plot -1.026 1.026   -1.399 .385 / 
\endpicture
\qquad\qquad
\beginpicture
\setcoordinatesystem units <1cm,1cm>         
\setplotarea x from -2 to 2, y from -2 to 2    
\put{$\bullet$} at .385 1.399
\put{$t_f$}[b] at .385 1.6
\put{$\bullet$} at -.385 1.399
\put{$s_1t_f$}[b] at -.385 1.6
\put{$\bullet$} at  -1.399 .385
\put{$s_1s_2s_1t_f$}[r] at  -1.6 .385
\put{$\bullet$} at -1.399 -.385
\put{$s_2s_1s_2s_1t_f$}[r] at -1.6 -.385
\put{$\bullet$} at 1.026 1.026
\put{$s_2t_f$}[bl] at 1.1 1.1
\put{$\bullet$} at -1.026 1.026
\put{$s_2s_1t_f$}[br] at -1.1 1.1
\plot .385 1.399   1.026 1.026 / 
\plot -1.026 1.026   -.385  1.399 / 
\plot -1.399 .385   -1.399 -.385 / 
\endpicture
$$
$$
\beginpicture
\setcoordinatesystem units <1cm,1cm>         
\setplotarea x from -2 to 2, y from -2 to 2    
\put{$\bullet$} at .385 1.399
\put{$t_g$}[b] at .385 1.6
\put{$\bullet$} at -.385 1.399
\put{$s_1t_g$}[b] at -.385 1.6
\put{$\bullet$} at -1.026 1.026
\put{$s_2s_1t_g$}[br] at -1.1 1.1
\put{$\bullet$} at  -1.399 .385
\put{$s_1s_2s_1t_g$}[r] at  -1.6 .385
\put{$\bullet$} at 1.026 1.026
\put{$s_2t_g$}[bl] at 1.1 1.1
\put{$\bullet$} at  1.399 .385
\put{$s_1s_2t_g$}[l] at  1.6 .385
\plot .385 1.399   1.026 1.026 / 
\plot 1.026 1.026  1.399 .385 / 
\plot -1.026 1.026   -.385  1.399 / 
\plot -1.026 1.026   -1.399 .385 / 
\endpicture
\qquad\qquad
\beginpicture
\setcoordinatesystem units <1cm,1cm>         
\setplotarea x from -2 to 2, y from -2 to 2    
\put{$\bullet$} at .385 1.399
\put{$t_h$}[b] at .385 1.6
\put{$\bullet$} at -.385 1.399
\put{$s_1t_h$}[b] at -.385 1.6
\put{$\bullet$} at -1.026 1.026
\put{$s_2s_1t_h$}[br] at -1.1 1.1
\put{$\bullet$} at  -1.399 .385
\put{$s_1s_2s_1t_h$}[r] at  -1.6 .385
\put{$\bullet$} at -1.399 -.385
\put{$s_2s_1s_2s_1t_h$}[r] at -1.6 -.385
\put{$\bullet$} at -1.026 -1.026
\put{$s_1s_2s_1s_2s_1t_h$}[tr] at -1.1 -1.1
\put{$\bullet$} at 1.026 1.026
\put{$s_2t_h$}[bl] at 1.1 1.1
\put{$\bullet$} at  1.399 .385
\put{$s_1s_2t_h$}[l] at  1.6 .385
\put{$\bullet$} at  1.399 -.385
\put{$s_2s_1s_2t_h$}[l] at  1.6 -.385
\put{$\bullet$} at 1.026 -1.026
\put{$s_1s_2s_1s_2t_h$}[tl] at 1.1 -1.1
\put{$\bullet$} at .385 -1.399
\put{$s_2s_1s_2s_1s_2t_h$}[tl] at .385 -1.6
\put{$\bullet$} at -.385 -1.399
\put{$s_1s_2s_1s_2s_1s_2t_h$}[tr] at -.385 -1.6
\plot .385 1.399   -.385 1.399 / 
\plot -1.026 1.026   -.385  1.399 / 
\plot -1.026 1.026   -1.399 .385 / 
\plot -1.399 -.385   -1.399 .385 / 
\plot -1.399 -.385   -1.026 -1.026 / 
\plot 1.026 1.026  1.399 .385 / 
\plot 1.399 -.385   1.399 .385  / 
\plot 1.399 -.385   1.026 -1.026  / 
\plot .385 -1.399   1.026 -1.026  / 
\plot .385 -1.399   -.385 -1.399  / 
\endpicture
$$
$$
\beginpicture
\setcoordinatesystem units <1cm,1cm>         
\setplotarea x from -2 to 2, y from -2 to 2    
\put{$\bullet$} at .385 1.399
\put{$t_i$}[b] at .385 1.6
\put{$\bullet$} at 1.026 1.026
\put{$s_2t_i$}[bl] at 1.1 1.1
\put{$\bullet$} at  1.399 .385
\put{$s_1s_2t_i$}[l] at  1.6 .385
\put{$\bullet$} at  1.399 -.385
\put{$s_2s_1s_2t_i$}[l] at  1.6 -.385
\put{$\bullet$} at 1.026 -1.026
\put{$s_1s_2s_1s_2t_i$}[tl] at 1.1 -1.1
\put{$\bullet$} at .385 -1.399
\put{$s_2s_1s_2s_1s_2t_i$}[tl] at .385 -1.6
\plot 1.026 1.026  1.399 .385 / 
\plot 1.399 -.385   1.399 .385  / 
\plot 1.399 -.385   1.026 -1.026  / 
\endpicture
\qquad\qquad
\beginpicture
\setcoordinatesystem units <1cm,1cm>         
\setplotarea x from -2 to 2, y from -2 to 2    
\put{$\bullet$} at .385 1.399
\put{$t_j$}[b] at .385 1.6
\put{$\bullet$} at -.385 1.399
\put{$s_1t_j$}[b] at -.385 1.6
\put{$\bullet$} at -1.026 1.026
\put{$s_2s_1t_j$}[br] at -1.1 1.1
\put{$\bullet$} at 1.026 1.026
\put{$s_2t_j$}[bl] at 1.1 1.1
\put{$\bullet$} at  1.399 .385
\put{$s_1s_2t_j$}[l] at  1.6 .385
\put{$\bullet$} at  1.399 -.385
\put{$s_2s_1s_2t_j$}[l] at  1.6 -.385
\plot .385 1.399   -.385 1.399 / 
\plot -1.026 1.026   -.385  1.399 / 
\plot 1.026 1.026  1.399 .385 / 
\plot 1.399 -.385   1.399 .385  / 
\endpicture
$$
\bigskip
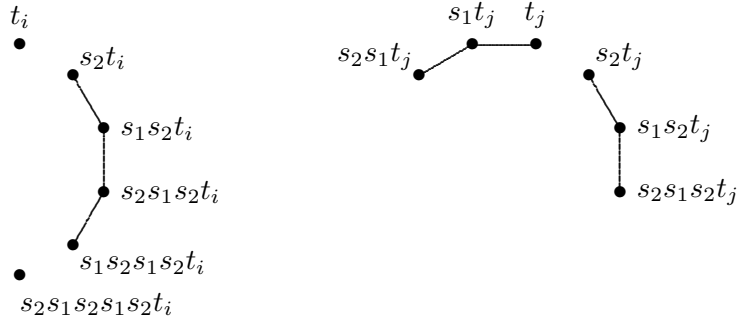
\centerline{{\bf Figure 6.2.} Calibration graphs for central characters in Table
6.1}

\subsection{The analysis.}

\smallskip\noindent
{\it Central characters $t_a$, $t_b$, $t_c$, $t_d$, $t_h$:}
The central characters $t_e$, $t_f$, $t_g$, $t_i$ and $t_j$ are the only ones which
have both $Z(t)$ and $P(t)$ nonempty.  The other cases are handled by
Theorem 1.16 and Theorem 1.11 as in the cases of central characters $t_a$,
$t_b$, $t_g$ and $t_o$ for type $A_2$.  

The Langlands parameters for each module can be determined from its weight
structure.  The indexing triple is determined from the Langlands data by using the
induction theorem of Kazhdan and Lusztig (see the discussion in [BM, p.34]).  Let
us give an example to illustrate the procedure.  The Langlands parameters
$(s_1s_2t_c,\{1\})$ for the $4$ dimensional representation
with central character $t_c$ correspond to the indexing triple
$(s_2t_c,e_{\alpha_1},1)$ which is conjugate to the
triple $(t_c,s_2e_{\alpha_1},1)=(t_c,e_{\alpha_1+3\alpha_2},1)$. 

The indexing triples for the tempered representations cannot be determined with
the use of the Kazhdan-Lusztig induction theorem.  The indexing triples
for the tempered representations with real central character are determined
from the Springer correspondence, see Table 6.3 and [BM, p.34].  The two tempered
representations with central characters $t_c$ and $t_d$ do not have
real central character. By the last two sentences of [Lu3, 2.10]
we know that the indexing triples for these representations contain
the subregular nilpotent and that the component groups are
isomorphic to $\ZZ/3\ZZ$ and $\ZZ/2\ZZ$ respectively.   In both cases the
component group acts trivially on $K(\cB_{s,u})$ and so $\rho=1$.  
The fact that the elements
$e_{\alpha_1}+e_{\alpha_1+3\alpha_2}$ and $e_{\alpha_1}+e_{\alpha_1+2\alpha_2}$
are representatives of the subregular nilpotent orbit can be derived from the
analysis in [Ja, Theorem 4.40] or [Sh].  This determines the triples
$(t_c, e_{\alpha_1}+e_{\alpha_1+3\alpha_2},1)$ and $(t_d,
e_{\alpha_1}+e_{\alpha_1+2\alpha_2},1)$.  

\medskip\noindent
{\it Central character $t_e$:}  
Theorem 1.11 applied to the placed skew shapes $(t_e, \{\alpha_1,
\alpha_1+\alpha_2\})$, $(t_e, \{\alpha_1\})$ and $(t_e, \{\alpha_1,
\alpha_1+\alpha_2,\alpha_1+2\alpha_2\})$ produces,
respectively, a two dimensional irreducible module $M$ with
$\supp(M)=\{s_2s_1t_e,s_1s_2s_1t_e\}$, a one dimensional irreducible module $N$
with
$\supp(N)=\{s_1t_e\}$ and a one dimensional irreducible module $N^*$ with
$\supp(N^*)=\{s_2s_1s_2s_1t_e\}$. 
Lusztig [Lu3] Theorem 4.7 constructs a 3-dimensional irreducible
$\tilde H$-module $P$ with $\dim(P_{t_e}^{\rm gen})=2$ and
$\dim(P_{s_1t_e}^{\rm gen})=1$.  In Lusztig's notation this is
the module labeled by the graph $\cG$ for $\tilde G_2$.  

As described in [Lu3, 4.23] we can twist the module $P$ by
an involutive automorphism of 
$\tilde H$ to obtain another 3-dimensional irreducible module
$P^*$ which has 
$\dim((P^*)_{s_2s_1s_2s_1t_e}^{\rm gen})=2$ and 
$\dim((P^*)_{s_1s_2s_1s_2s_1t_e}^{\rm gen})=1$.

All of the modules $M$, $N$, $P$, $N^*$, $P^*$ must appear
as composition factors of the principal series module $M(t_e)$.
By comparing dimensions of weight spaces, any other module $Q$
which appears in a composition series of $M(t_e)$ must have
$\supp(Q)\subseteq \{s_2s_1t_e,s_1s_2s_1t_e\}$.  Theorem 1.11(b) then 
implies that $Q$ must be isomorphic to $M$. Thus Theorem 1.15 implies that
$M$, $N$, $P$, $N^*$, and $P^*$ are (up to isomorphism) all the irreducible modules
with central character $t_e$.

The Langlands parameters for each module are determined from its weight
structure.  The Kazhdan-Lusztig induction theorem allows us to use the
Langlands parameters to determine the
indexing triples for the modules which are not tempered.  Since $t_e$ is a
real central character the indexing triples for the tempered representations
can be determined from the Springer correspondence, see Table 6.3.
Alternatively, one can get these triples from [Lu3, 2.10] where it is explained
that the nilpotent in the indexing triple is subregular, the variety $\cB_{s,u}$
(where $s=t_e$ and $u$ is subregular) consists of three disjoint points and a
projective line, and the component group is isomorphic to the symmetric group
$S_3$.  The symmetric group
$S_3$ acts trivially on the line and permutes the three points, which implies that
the line corresponds to
$\rho=(3)$ (trivial representation of $S_3$) and the three points are split 
between the isotypic components $\rho=(3)$ and $\rho=(21)$.  In this case
the standard module $M_{s,u,(1^3)}=0$.  The projective line in $\cB_{s,u}$
corresponds to the two dimensional weight space $(P^*)_{t_e}^{\rm gen}$ in the
module $P^*$.

\smallskip\noindent
{\it Central character $t_f$:}  
Let $\CC v_f$ and $\CC v_{s_1f}$ be the one dimensional
representations of $\tilde H_{\{1\}}$ given by
$$
\matrix{
T_1v_f=qv_f,  &&X^{\alpha_1}v_f=q^2v_f,
&&X^{3\alpha_2}v_f=q^{-2}v_f, \cr 
T_1v_{s_1f}=-q^{-1}v_{s_1f}, &&X^{\alpha_1}v_{s_1f}=q^{-2}v_{s_1f},
&&X^{3\alpha_2}v_{s_1f}=q^{4}v_{s_1f}. \cr }
$$
Let 
$$M_1 = \Ind_{\tilde H_{\{1\}}}^{\tilde H}(\CC v_f)
\qquad\hbox{and}\qquad
M_2 = \Ind_{\tilde H_{\{1\}}}^{\tilde H}(\CC v_{s_1f}).$$
By Lemma 1.17 these modules have weights
$\supp(M_1)=\{s_2t_f,t_f,s_1t_f,s_2s_1t_f\}$ and
$\supp(M_2)=\{s_1t_f,s_2s_1t_f,$ $s_1s_2s_1t_f,s_2s_1s_2s_1t_f\}$ respectively. 
Both $M_1$ and $M_2$ are $6$ dimensional.

Let $M$ be any $\tilde H$-module such that $M_{s_2t_f}\ne 0$.
By Lemma 1.19 and Proposition 1.6, $\dim(M_{t_f}^{\rm gen})
=\dim(M_{s_2t_f}^{\rm gen})\ge 2$.  Since
$\langle s_2\alpha_1,\alpha_1^\vee\rangle =\langle
\alpha_1+3\alpha_2,\alpha_1^\vee\rangle \ne 0$ it follows from Lemma 1.19 (b) that
$\dim(M_{s_1t_f}^{\rm gen})\ge 1$.  Then, by Proposition 1.6, 
$\dim(M_{s_2s_1t_f}^{\rm gen})\ge 1$.  Adding these numbers up we see that
$\dim(M)\ge 6$.  It follows that $M_1$ is irreducible.
An analogous argument can be applied to conclude that $M_2$ is irreducible.

\smallskip\noindent
{\it Central character $t_g$:} 
Let $\CC v_g$ and $\CC v_{s_1g}$ be the one dimensional
representations of $\tilde H_{\{1\}}$ given by
$$
\matrix{
T_1v_g=qv_g,  &&X^{\alpha_1}v_g=q^2v_g,
&&X^{2\alpha_2}v_g=q^{-2}v_g, \cr 
T_1v_{s_1g}=-q^{-1}v_{s_1g}, &&X^{\alpha_1}v_{s_1g}=q^{-2}v_{s_1g},
&&X^{2\alpha_2}v_{s_1g}=q^2v_{s_1g}. \cr }
$$
Let 
$$M_1 = \Ind_{\tilde H_{\{1\}}}^{\tilde H}(\CC v_g)
\qquad\hbox{and}\qquad
M_2 = \Ind_{\tilde H_{\{1\}}}^{\tilde H}(\CC v_{s_1g}).$$
By Lemma 1.17 these modules have weights
$\supp(M_1)=\{s_1s_2t_g,s_2t_g,t_g\}$ and
$\supp(M_2)=\{s_1t_g,s_2s_1t_g,$ $s_1s_2s_1t_g\}$ respectively.  Both
$M_1$ and $M_2$ are $6$ dimensional.

Let $M$ be any $\tilde H$-module such that $M_{s_1s_2t_g}\ne 0$.
By Lemma 1.19 and Proposition 1.6,
$\dim(M_{t_g}^{\rm gen})=\dim(M_{s_2t_g}^{\rm gen})
=\dim(M_{s_1s_2t_g}^{\rm gen})\ge 2$. Thus
$\dim(M)\ge 6$.  It follows that $M_1$ is irreducible. An analogous argument can
be applied to conclude that $M_2$ is irreducible.

\smallskip\noindent
{\it Central character $t_i$:}  
Let $\CC v_i$ and $\CC v_{s_2i}$ be the one dimensional
representations of $\tilde H_{\{2\}}$ given by
$$
\matrix{
T_2v_i=qv_i,  &&X^{\alpha_1}v_i=v_i,
&&X^{\alpha_2}v_i=q^2v_i, \cr 
T_2v_{s_2i}=-q^{-1}v_{s_2i}, &&X^{\alpha_1}v_{s_2i}=q^6v_{s_2i},
&&X^{\alpha_2}v_{s_2i}=q^{-2}v_{s_2i}. \cr }
$$
Let 
$$M_1 = \Ind_{\tilde H_{\{2\}}}^{\tilde H}(\CC v_i)
\qquad\hbox{and}\qquad
M_2 = \Ind_{\tilde H_{\{2\}}}^{\tilde H}(\CC v_{s_2i}).$$
An argument similar to that for the central character $t_f$ 
shows that $M_1$ and $M_2$ are irreducible.

\smallskip\noindent
{\it Central character $t_j$:}  
Let $\CC v_j$ and $\CC v_{s_2j}$ be the one dimensional
representations of $\tilde H_{\{2\}}$ given by
$$
\matrix{
T_2v_j=qv_j,  &&X^{2\alpha_1}v_j=q^{-6}v_j,
&&X^{\alpha_2}v_j=q^2v_j, \cr 
T_2v_{s_2j}=-q^{-1}v_{s_2j}, &&X^{2\alpha_1}v_{s_2j}=q^6v_{s_2j},
&&X^{\alpha_2}v_{s_2j}=q^{-2}v_{s_2j}. \cr }
$$
Let 
$$M_1 = \Ind_{\tilde H_{\{2\}}}^{\tilde H}(\CC v_j)
\qquad\hbox{and}\qquad
M_2 = \Ind_{\tilde H_{\{2\}}}^{\tilde H}(\CC v_{s_2j}).$$
An argument similar to that for the central character $t_g$ 
shows that $M_1$ and $M_2$ are irreducible.

\vfill\eject

\centerline{\smallcaps References}

\medskip


\medskip
\item{[BM]} {\smallcaps D. Barbasch and A. Moy}, {\it A unitarity criterion for
$p$-adic groups}, Invent. Math. {\bf 98} (1989), 19-37.

\medskip
\item{[Bou]} {\smallcaps N. Bourbaki}, {\it Groupes et alg\`ebres de Lie,
Chapitres 4,5 et 6}, Elements de Math\'ematique, Hermann, Paris 1968. 

\medskip
\item{[Ca]} {\smallcaps R.W.\ Carter}, 
{\it Finite Groups of Lie Type: Conjugacy Classes and Complex Characters}, 
John Wiley and Sons, 1985.

\medskip
\item{[CG]} {\smallcaps N.\ Chriss and V.\ Ginzburg}, 
{\sl Representation Theory and Complex Geometry}, Birkh\"auser, 1997.

\medskip
\item{[Ev]} {\smallcaps S. Evens}, {\it The Langlands classification for graded
Hecke algebras}, Proc. Amer. Math. Soc. {\bf 124} (1996), 1285--1290.

\medskip
\item{[HO1]} {\smallcaps G.J.\ Heckman and E.M.\ Opdam}, 
{\it Yang's system of particles and Hecke algebras}, 
Ann. of Math. (2) {\bf 145} (1997), 139--173.

\medskip
\item{[HO2]} {\smallcaps G.J.\ Heckman and E.M.\ Opdam}, 
{\it Harmonic analysis for affine Hecke algebras}, in Current Developments
in Mathematics, Intern. Press, Boston, 1996.

\medskip
\item{[Ja]} {\smallcaps D.\ Jackson},  {\sl On nilpotent orbits of type $G_2$},
Ph. D. Thesis, University of Sydney, 1997.

\medskip
\item{[Ka]} {\smallcaps S-i.\ Kato},
{\it Irreducibility of principal series representations for Hecke
algebras of affine type}, 
J. Fac. Sci. Univ. Tokyo Sect. IA Math. {\bf 28} (1981), 929--943.

\medskip
\item{[KL]} {\smallcaps D.\ Kazhdan and G.\ Lusztig}, 
{\it Proof of the Deligne-Langlands conjecture for Hecke algebras}, 
Invent. Math. {\bf 87} (1987), 153--215.

\medskip
\item{[Kr]} {\smallcaps C. Kriloff}, {\it  Some interesting nonspherical tempered
representations of graded Hecke algebras}, Trans. Amer. Math. Soc., to appear.

\medskip
\item{[Lu1]} {\smallcaps G.\ Lusztig},
{\it Singularities, character formulas, and a $q$-analog of weight
multiplicities}, Analysis and topology on singular spaces, II, III (Luminy, 1981),
 Ast\'erisque {\bf 101-102}, Soc. Math. France, Paris, 1983, 208--229.

\medskip
\item{[Lu2]} {\smallcaps G.\ Lusztig},
{\it Affine Hecke algebras and their graded version}, J. Amer. Math.
Soc. {\bf 2} (1989), 599--635.

\medskip
\item{[Lu3]} {\smallcaps G.\ Lusztig}, 
{\it Some examples of square integrable
representations of semisimple $p$-adic groups}, 
Trans.\ Amer.\ Math.\ Soc.\ {\bf 277} (1983), 623--653.

\medskip
\item{[Lu4]} {\smallcaps G. Lusztig},
{\it Equivariant K-theory and representations of Hecke algebras}, Proc. Amer.
Math. Soc. {\bf 94} (1985), 337--342.


\medskip
\item{[Ma]} {\smallcaps H.\ Matsumoto}, 
{\sl Analyse harmonique dans les syst\`emes de Tits bornologiques de type affine},
Lect. Notes in Math. {\bf 590}, Springer-Verlag, Berlin-New York,
1977.

\medskip
\item{[Ra1]} {\smallcaps A. \ Ram}, 
{\it Calibrated representations of affine Hecke algebras}, preprint 1998.

\medskip
\item{[Ra2]} {\smallcaps A.\ Ram}, 
{\it Standard Young tableaux for finite root systems}, preprint 1998.

\medskip
\item{[Ra3]} {\smallcaps A.\ Ram}, 
{\it Skew shape representations are irreducible}, in preparation.

\medskip
\item{[RR1]} {\smallcaps A.\ Ram and J. Ramagge}, 
{\it  Jucys-Murphy elements come from affine Hecke algebras}, in preparation.

\medskip
\item{[RR2]} {\smallcaps A.\ Ram and J. Ramagge}, {\it Calibrated
representations and the $q$-Springer correspondence}, in preparation.

\medskip
\item{[Ro]} {\smallcaps F.\ Rodier}, 
{\it D\'ecomposition de la s\'erie principale des groupes r\'eductifs $p$-adiques},
Noncommutative harmonic analysis and Lie groups (Marseille, 1980),
Lect. Notes in Math. {\bf 880}, Springer, Berlin-New York, 1981, 408--424.

\medskip
\item{[St]} {\smallcaps R.\ Steinberg}, 
{\it Endomorphisms of linear algebraic groups}, Mem. Amer.
Math. Soc. {\bf 80}, Amer. Math. Soc., Providence, R.I. 1968.

\medskip
\item{[Sh]} {\smallcaps U. Stuhler}, {\it Unipotente und nilpotente
Klassen in einfachen Gruppen und Liealgebren vom typ $G_2$}, Indag.
Math. {\bf 74} (1971), 365--378.

\medskip
\item{[Xi]} {\smallcaps N.\ Xi}, 
{\it  Representations of affine Hecke algebras}, 
Lect. Notes in Math. {\bf 1587}, Springer-Verlag, Berlin, 1994.

\vfill\eject
\end

\subsection{Restriction to $\tilde H(Q)$}

Every weight $t\in T(Q)$, $t(X^{\alpha_1})=t_1^2$ has two lifts in $T(P)$,
$$t(X^{\omega_1})=\pm t_1.$$

The two irreducible representations corresponding to the central
characters $t(X^{\omega_1})=\pm 1$ restrict to the same
representation of $\tilde H(Q)$.

The weight $t$ given by $t(X^{\omega_1})=i$ is regular and 
the principal series representation $M(t)$ is irreducible,
$2$-dimensional, and calibrated.
Upon restriction to $\tilde H(Q)$, the module $M(t)$
splits into two one-dimensional calibrated representations with central character
$t$ given by $t(X^{\alpha_1})=-1$.  Note that
$(s_1t)(X^{\omega_1})=-i$ which is the other square root of
$t(X^{\alpha_1})=-1$.

\subsection{The KL labeling}

Let $G=SL_2(\CC)$ and choose $B$ and $T$ to be the subgroups
of upper triangular and diagonal matrices, respectively.
The weight given by $t(X^{\omega_1}) = t_1$ is identified with the matrix 
$$t = 
\pmatrix{t_1 &0 \cr 0 &t_1^{-1} \cr}$$
and $u$ is the regular unipotent 
$$u = \pmatrix{1 &1\cr 0 &1 \cr}.$$

$$\eqalign{
X^{\omega_1} v_{\pm q} &= \pm q v_{\pm q}, \cr
T_1          v_{\pm q} &= q v_{\pm q}, \cr
}
$$
$$\eqalign{
X^{\omega_1} v_{\pm q^{-1}} &= \pm q^{-1} v_{\pm q^{-1}}, \cr
T_1          v_{\pm q^{-1}} &= -q^{-1} v_{\pm q^{-1}}. \cr
}
$$